\documentclass{amsart}
\usepackage{amsmath,amssymb,amsthm}
\usepackage{accents}
\usepackage{color}
\usepackage{mathtools}
\mathtoolsset{showonlyrefs}
\usepackage[unicode]{hyperref}
\hypersetup{pdftitle=viscoelastic fluids with stress diffusion,
            pdfauthor=Michal Bathory,
            hidelinks,
            pdfstartview=FitH
}

\newtheorem*{theorem*}{Theorem}

\theoremstyle{definition}
\newtheorem*{definition*}{Definition}

\newtheorem*{remark*}{Remark}

\newcommand{\eps}{\varepsilon}

\newcommand{\R}{\mathbb{R}}

\newcommand{\N}{\mathbb{N}}

\newcommand{\f}[2]{\frac{#1}{#2}}
\newcommand{\tf}[2]{\tfrac{#1}{#2}}

\newcommand{\dd}[1]{\,\mathrm{d}{#1}}

\newcommand{\embl}{\hookrightarrow}

\DeclareMathOperator{\Div}{div}

\DeclareMathOperator{\tr}{tr}

\DeclareMathOperator{\sym}{sym}

\DeclareMathOperator{\spa}{span}

\DeclareMathOperator*{\esssup}{ess\ sup}

\def\XXint#1#2#3{{\setbox0=\hbox{$#1{#2#3}{\int}$}
\vcenter{\hbox{$#2#3$}}\kern-.5\wd0}}

\usepackage{mleftright,xparse}

\NewDocumentCommand\xDeclarePairedDelimiter{mmm}
 {%
  \NewDocumentCommand#1{som}{%
   \IfNoValueTF{##2}
    {\IfBooleanTF{##1}{#2##3#3}{\mleft#2##3\mright#3}}
    {\mathopen{##2#2}##3\mathclose{##2#3}}%
  }%
 }
\xDeclarePairedDelimiter{\abs}{\lvert}{\rvert}
\xDeclarePairedDelimiter{\set}{\lbrace}{\rbrace}
\xDeclarePairedDelimiter{\norm}{\lVert}{\rVert}
\xDeclarePairedDelimiter{\scal}{\langle}{\rangle}

\newcommand{\normm}[1]{{\left\vert\kern-0.25ex\left\vert\kern-0.25ex\left\vert #1
		\right\vert\kern-0.25ex\right\vert\kern-0.25ex\right\vert}}

\newcommand{\bs}{\boldsymbol}
\newcommand{\bb}{\mathbb}
\newcommand{\D}{\bb D}

\newcommand{\Wb}{\bb W}

\newcommand{\JJ}{\bb J}
\newcommand{\ve}{\bs v}

\newcommand{\we}{\bs w}

\newcommand{\n}{\bs n}

\newcommand{\fe}{\bs f}

\newcommand{\fit}{\bs \varphi}

\newcommand{\ii}{\int_{\Omega}}

\newcommand{\nn}{\nabla}

\newcommand{\LH}{L^2_{\n,\Div}}

\newcommand{\T}{\mathbb T}

\newcommand{\z}{\bs z}

\newcommand{\I}{\mathbb I}
\newcommand{\Dv}{\D\ve}

\newcommand{\wc}{\rightharpoonup}
\newcommand{\wcs}{\overset{*}{\rightharpoonup}}

\newcommand{\ot}{\otimes}

\newcommand{\la}{\lambda}
\newcommand{\La}{\Lambda}

\newcommand{\B}{\mathbb{B}}

\newcommand{\je}{\bs j}
\newcommand{\der}[1]{\accentset{\diamond}{#1}}
\newcommand{\A}{\mathbb{A}}

\renewcommand{\dot}[1]{\accentset{\bullet}{#1}}
\newcommand{\pp}{\partial}

\newcommand{\Log}{\bb{H}}

\newcommand{\WH}{W^{1,2}}
\newcommand{\WHd}{W^{-1,2}}
\newcommand{\Wn}[1]{W_{\n}^{1,#1}}
\newcommand{\Wnd}[1]{W_{\n,\Div}^{1,#1}}
\newcommand{\Wndt}{W_{\n,\Div}^{3,2}}
\newcommand{\Wndd}[1]{W_{\n,\Div}^{-1,#1}}
\newcommand{\Wnddt}{W_{\n,\Div}^{-3,2}}
\newcommand{\PD}{\R^{3\times 3}_{>0}}
\newcommand{\Tr}{\mathcal{T}}
\newcommand{\Sb}{\bb S}

\title[Mathematical analysis of Navier-Stokes-Oldroyd-B fluids]{Large data existence theory for three-dimensional unsteady flows of rate-type viscoelastic fluids with stress diffusion}

\author[M. Bathory]{Michal Bathory}
\address{Mathematical Institute of Charles University\\
	Faculty of Mathematics and Physics\\
	Charles University\\
	Sokolovsk\'a 83\\186 75 Praha 8\\Czech Republic}
\email{bathory@karlin.mff.cuni.cz}
\author[M. Bul\'i\v cek]{Miroslav Bul\'i\v cek}
\address{Mathematical Institute of Charles University\\
	Faculty of Mathematics and Physics\\
	Charles University\\
	Sokolovsk\'a 83\\186 75 Praha 8\\Czech Republic}
\email{mbul8060@karlin.mff.cuni.cz}
\author[J. M\'{a}lek]{Josef M\'{a}lek}
\address{Mathematical Institute of Charles University\\
	Faculty of Mathematics and Physics\\
	Charles University\\
	Sokolovsk\'a 83\\186 75 Praha 8\\Czech Republic}
\email{malek@karlin.mff.cuni.cz}

\thanks{Michal Bathory has been supported by Charles University Research program UNCE/SCI/023 and also by the project No.\ 1652119 financed by the Charles University Grant Agency (GAUK). Miroslav Bul\'{\i}\v{c}ek and Josef M\'{a}lek acknowledge the support of the project  No.\ 18-12719S financed by Czech science foundation (GA\v{C}R). Miroslav Bul\'{\i}\v{c}ek and Josef M\'{a}lek are members of the Jind\v{r}ich Ne\v{c}as center for mathematical modelling}
\subjclass{Primary 35Q35, 76A05, 76A10}
\keywords{viscoelasticity; viscoleastic fluid; Oldroyd-B; Johnson-Segalman; existence; weak solution; stress diffusion}

\numberwithin{equation}{section}

\begin{document}

\begin{abstract}
We prove that there exists a weak solution to a system governing an unsteady flow of a viscoelastic fluid in three dimensions, for arbitrarily large time interval and data. The fluid is described by the incompressible Navier-Stokes equations for the velocity $\ve$, coupled with a diffusive variant of a combination of the Oldroyd-B and the Giesekus models for a tensor $\B$. By a proper choice of the constitutive relations for the Helmholtz free energy (which, however, is non-standard in the current literature, despite the fact that this choice is well motivated from the point of view of physics) and for the energy dissipation, we are able to prove that $\B$ enjoys the same regularity as $\ve$ in the classical three-dimensional Navier-Stokes equations. This enables us to handle any kind of objective derivative of $\B$, thus obtaining existence results for the class of diffusive Johnson-Segalman models as well. Moreover, using a suitable approximation scheme, we are able to show that $\B$ remains positive definite if the initial datum was a positive definite matrix (in a pointwise sense). We also show how the model we are considering can be derived from basic balance equations and thermodynamical principles in a natural way.
\end{abstract}

	\maketitle
	
\section{Introduction}
	We aim to establish a global-in-time and  large-data existence theory, within the context of weak solutions, to a class of homogeneous incompressible rate-type viscoelastic fluids flowing in a closed three-dimensional container. The studied class of models can be seen as the Navier-Stokes system (for which a similar existence theory is well known, cf.~\cite{Leray}) coupled with a viscoelastic rate-type fluid model that shares the properties of both Oldroyd-B and Giesekus models and is completed with a diffusion term. Such models are frequently encountered in the theory of non-Newtonian fluid mechanics, see \cite{Elkareh,Dostalik2019} and further references cited in \cite{Dostalik2019}.
	
In order to precisely formulate the problems investigated in this study, we start by introducing the necessary notation. For a bounded domain $\Omega\subset\R^3$ with the Lipschitz boundary $\partial \Omega$ and a time interval of the length $T>0$, we define the time-space cylinder $Q:=(0,T)\times\Omega$ and we also set $\Sigma:=(0,T) \times \partial \Omega$ for a part of its boundary. The symbol $\n$ denotes the outward unit normal vector on $\partial\Omega$ and, for any vector~$\z$, the vector~$\z_\tau$ denotes the projection of the vector to a tangent plane on $\partial \Omega$, i.e., $\z_{\tau}:= \z - (\z \cdot \n)\n$. Then, for a given density of the external body forces $\fe:Q\to\R^3$, a given initial velocity $\ve_0:\Omega\to\R^3$ and a given initial extra stress tensor $\B_0:\Omega\to\PD$ (here $\PD$ denotes the set of symmetric positive definite ($3\times 3$)-matrices), we look for a vector field $\ve:Q\to\R^3$, a scalar field $p:Q\to\R$ and a positive definite matrix field $\B:Q\to\PD$ solving the following system in~$Q$:
\begin{align}
	\Div\ve&=0,\label{divv}\\
	\pp_t\ve+(\ve\cdot\nn)\ve-\nu\Delta\ve+\nn p&=2\mu a\Div({(1\!-\!\gamma)}(\B-\I)+\gamma(\B^2\!-\!\B))+\fe,\label{vv}
\end{align}
\begin{equation}
\begin{split}
	\pp_t\B+(\ve\cdot\nn)\B&+\delta_1(\B-\I)+\delta_2(\B^2-\B)-\la\Delta\B\\
	&=\tf{a+1}2(\nn\ve\B+(\nn\ve\B)^T)+\tf{a-1}2(\B\nn\ve+(\B\nn\ve)^T),\label{BB}
\end{split}
\end{equation}
and being completed by the following boundary conditions on $\Sigma$:
\begin{equation}
\begin{split}
	\ve\cdot\n&=0,\\
	-\sigma\ve_{\tau}&=\left(\left(\nu\nn\ve+\nu(\nn\ve)^T+2\mu a(1\!-\!\gamma)(\B-\I)+2\mu a\gamma(\B^2-\B)\right)\n\right)_{\tau},\\
	(\n\cdot\nn)\B&=\mathbb{O},  \qquad\qquad (\textrm{here } \mathbb{O} \textrm{ stands  for zero } 3\times 3\textrm{-matrix}) \label{noflux}
	\end{split}
\end{equation} 
and by the initial conditions in $\Omega$:
	\begin{align}
	\ve(0,\cdot)&=\ve_0,\label{ic1}\\
	\B(0,\cdot)&=\B_0\label{ic2}.
	\end{align}
The parameters  $\gamma\in(0,1)$, $\nu, \la,\sigma>0$, $\delta_1,\delta_2\geq0$ and $a\in\R$ are given numbers. 

The main result of this study can be stated as:

{\it
Let $\ve_0$ and $\B_0$ be such that the initial total energy is bounded. Then, for sufficiently regular $\fe$, there exists a global-in-time weak solution to \eqref{divv}--\eqref{ic2}.
}

Although the above result is stated vaguely, we would like to emphasize that we are going to establish the {\bf long-time} existence of a weak solution for {\bf large data} and for {\bf three-dimensional} flows. A more precise and rigorous version of the above result including the correct function spaces and the properly defined weak formulation is stated in the Theorem below, see Section~\ref{secmain}.

We complete the introductory part by providing the physical background relevant to the studied problem and by recalling earlier results relevant to the problem \eqref{divv}--\eqref{ic2} analyzed here.

\subsection{Mathematical and physical background}

The system \eqref{divv}--\eqref{noflux} can be rewritten into a more concise form once one recognizes some physical quantities. First of all, let
	$$\Dv=\tf12(\nn\ve+(\nn\ve)^T)\quad\text{and}\quad\Wb\ve=\tf12(\nn\ve-(\nn\ve)^T)$$
	denote the symmetric and antisymmetric parts of the velocity gradient $\nn\ve$, respectively. Then, looking at the equation \eqref{vv},
	we see that \eqref{vv} is obtained from a general form of the balance of linear momentum, namely
	\begin{equation}
	\varrho \dot\ve = \Div\T + \varrho \fe, \label{p_BLM}
	\end{equation}
	once we set the density $\varrho=1$ and require that the  Cauchy stress tensor $\T$ has the form
	\begin{equation}\label{T}
	\T=-p\I+2\nu\Dv+2a\mu((1-\gamma)(\B-\I)+\gamma(\B^2-\B)).
	\end{equation}
	In \eqref{p_BLM}, $\dot\ve$ stands for the material time derivative of $\ve$, i.e., $\dot\ve=\pp_t\ve+(\ve\cdot\nn)\ve$.
Defining similarly the material time derivative of a tensor $\B$ as
$$
\dot\B=\pp_t\B+(\ve\cdot\nn)\B,
$$
we can recognize the presence of a general objective derivative in \eqref{BB}. Namely, defining
$$
\der{\B}=\dot{\B}-a(\Dv\B+\B\Dv)-(\Wb\ve\B-\B\Wb\ve),
$$
we can rewrite the system \eqref{divv}--\eqref{BB} into a more familiar form  as
\begin{align}
	\Div\ve&=0,\label{r1}\\
	\dot{\ve}&=\Div\T+\fe,\label{r2}\\
	\der{\B}+\delta_1(\B-\I)+\delta_2(\B^2-\B)&=\la\Delta\B,\label{r3}
\end{align}
which is supposed to hold true in $Q$ and which is completed by the initial conditions \eqref{ic1}, \eqref{ic2} fulfilled in $\Omega$ and by the boundary conditions \eqref{noflux} on $\Sigma$ that take the form:  
\begin{align}
	\ve\cdot\n&=0,\label{r4}\\
	(\T\n)_{\tau}&=-\sigma\ve_{\tau},\label{r5}\\
	(\n\cdot\nn)\B&=\mathbb{O}.\label{r6}
\end{align}

We provide several comments regarding \eqref{T}--\eqref{r3} as well as the boundary conditions \eqref{r4}--\eqref{r6}. The Navier slip boundary condition \eqref{r5} (and in general
all boundary conditions allowing the
fluid to slip ever so slightly) has recently attracted lot of attention. It was well documented that in certain situations the Navier slip boundary conditions are more appropriate than no slip boundary conditions, we refer e.g.\ to \cite{Drda,Hatzi,Jager,stick1} or \cite{Zaba4} and references therein. In addition, it was shown that the Navier slip boundary condition can be understood as an asymptotic limit of no slip boundary conditions in case we consider rough and highly oscillating boundary, see e.g.~\cite{Amirat,Basson,Bucur}. Furthermore, for the classical Navier-Stokes equation or the Stokes equation, we can say that the available mathematical theory for no-slip boundary condition has been already ``re-proven" for Navier boundary conditions, see e.g.~\cite{Amrouche2011} for the existence analysis, \cite{Amrouche2014,Baba2019,Macha2014} for regularity theory for the Stokes system and \cite{Beirao2020} for a conditional regularity result for Navier-Stokes system. The key difference and also the main mathematical advantage of the Navier slip boundary conditions is, that for smooth domains, namely if $\Omega\in \mathcal{C}^{1,1}$, we can introduce the pressure $p$ as an integrable function, e.g., by using an additional layer of approximation as in \cite{Bulicek2009}, see also \cite{Zaba2,Zaba3} or \cite{Blechta2019} which discuss the treatment of the pressure in evolutionary models subject to the Navier boundary condition. Nevertheless, since we shall always deal with formulation without the pressure (see the Definition), we can also treat the Dirichlet boundary condition, as well as very general implicitly specified boundary conditions see e.g.~\cite{Zaba4,stick1,stick2} or \cite{Blechta2019}. The Neumann boundary condition for $\B$ is considered here only for simplicity and without any specific physical meaning.
	
A further aspect, which makes the above system more complicated than the Navier-Stokes equation is the form of the Cauchy stress tensor $\T$ as in \eqref{T}. The term $-p\I+2\nu\Dv$ corresponds to the standard Newtonian fluid flow model with a constant kinematic viscosity $\nu$. The next part of the Cauchy stress, which depends linearly on $\B$, appears in all the viscoelastic rate-type fluid models - see, e.g., \cite[(7.20b), (8.20e)]{Malek2018}, \cite[(6.43e)]{Hron2017} or \cite[(43a)]{Dostalik2019}. On the other hand, the addition of the term $2a\mu \gamma(\B^2-\B)$ is, to our best knowledge, considered here for the first time. The fact that we require that $\gamma$ is positive (and strictly less than $1$) plays a key role in the analysis of the problem, as will be shown below. Note that the linearization of $\T$ with respect to $\B$ when $\B$ is close to the identity $\I$ yields $$\T=-p\I+2\nu\Dv+2a\mu(\B-\I)$$
and we recover the standard form of $\T$ (after possible redefinition of the pressure).
	
The quantity $\B$ 
takes into account the elastic responses of the fluid and  the equation \eqref{r3} describes its evolution in the current configuration (Eulerian coordinates), just as the velocity
$\ve$. It is frequent to call the tensor $\mu(\B-\I)$ the extra stress or conformation tensor and to denote it by $\bs{\tau}$. More importantly, since the material derivative of $\B$ is not objective, it must be ``corrected" and this is the reason, why in \eqref{r3} the derivative $\der{\B}$ appears. The parameter $a$ in the definition of $\der{\B}$ determines the type of the objective derivative. The case $a=1$ leads to the upper convected Oldroyd derivative, that has favourable physical properties and that leads to a clear interpretation of $\B$ within the thermodynamical framework developed in \cite{RS_2000}, see also \cite{RS_2004, MaRaTu15, MaRaTu18, MaPr18}. Next, the case $a=0$ leads to the corrotational Jaumann-Zaremba derivative and this is the only case for which the analysis is much simpler than in other cases. Furthermore, if $a\in[-1,1]$, one obtains the entire class of Gordon-Schowalter derivatives. However, it turns out that the physical properties of these derivatives are irrelevant for the analysis presented below (except the case $a=0$), therefore we may take any $a\in\R$. For $a=1$ and $\la=0$ we distinguish two cases: if $\delta_1>0$ and $\delta_2=0$ we obtain the classical Oldroyd-B model while if $\delta_1=0$ and $\delta_2>0$ we get the Giesekus model. Next, by considering $a\in[-1,1]$, we obtain the class of Johnson-Segalman models. If we further let $\la>0$, we are introducing diffusive variants of the previous models. It has been observed that including the diffusion term in \eqref{r3} is physically reasonable, see, e.g., \cite{Elkareh} or \cite{Dostalik2019} and references therein. However, up to now, it has been unknown what precise form should the diffusion term take and also whether it actually helps in the analysis of the model. Our main result provides a partial answer to this question, namely: for $\gamma\in(0,1)$ and with the diffusion term being of  the form $\Delta \B$ (or more generally, a linear second order operator), the global existence of a weak solution is available.

The reader familiar with the equations describing flows of the standard \mbox{Oldroyd-B} viscoleastic rate-type fluid can identify two deviations in the set of equations \eqref{r1}--\eqref{r3} studied hereafter. We provide a few comments on these differences.

The first deviation concerns the incorporation of the stress diffusion term, i.e.\ the term $-\Delta\B$, into the equations. Following the pioneering work of \cite{Elkareh} it is clear that a quantity related to $|\nn\B|^2$ has to be added into the list of underlying dissipation mechanisms. On the other hand, the precise form in which stress diffusion should appear depends on the choice of a thermodynamical approach and specific assumptions. In fact, using the thermodynamical concepts as in \cite{Malek2018} or \cite{Dostalik2019}, one can derive models, where the stress diffusion term takes the form $-\B\Delta\B-\Delta\B\B$, $-\B^{\f12}\Delta\B\B^{\f12}$ etc., however, we would prefer $-\Delta\B$ simply because it coincides with the form proposed by \cite{Elkareh}, and, from the perspective of PDE analysis and numerical approximation, one prefers to deal with stress diffusion that leads to a linear operator.

The second deviation from usual viscoelastic models consists in the presence of the term $(\B^2-\B)$ in the Cauchy stress tensor, see \eqref{T}. This term arises if we slightly modify energy storage mechanism and apply the thermodynamic approach as developed in \cite{Malek2018}. In what follows, we shall give a clear interpretation and a thermodynamic derivation of our model.

\subsection{Thermodynamical derivation of the model}\label{shoowed}
Viscoelastic models with (nonlinear) stress diffusion, but without the term $\B^2$ in the stress tensor are derived, e.g., in \cite{Malek2018} and \cite{Dostalik2019} even in the temperature-dependent case. Here, we will briefly explain the approach in a simplified isothermal setting (sufficient for the purpose of this study), referring to the cited works for the derivation in a complete thermal setting and for more details.

First, we postulate the constitutive equation for the Helmholtz free energy in the form
\begin{equation}\label{helm}
\psi(\B):=\mu((1-\gamma)(\tr\B-3-\ln\det\B)+\tf12\gamma|\B-\I|^2),
\end{equation}
where $\mu >0$ and $\gamma\in [0,1]$ is a parameter interpolating between two forms of the energy.
The choice $\gamma=0$ would lead to a standard Oldroyd-B diffusive model. To our best knowledge, the case $\gamma>0$ was not considered before in literature. The term $\tf12\gamma|\B-\I|^2$, which is newly included in $\psi$ is obviously convex with the minimum at $\B=\I$ and depends only on $\tr \B$ and on $\tr (\B \B)$, i.e., on invariants of $\B$, therefore it does not violate any of the basic principles of continuum physics. Moreover, such an addition does not affect the first three terms in the asymptotic expansion of $\psi$ near $\I$, on the logarithmic scale. To see this, let $\Log$ denote the Hencky logarithmic tensor satisfying $e^{\Log}=\B$ (which exists due to the positive definiteness of $\B$). Using Jacobi's identity, we compute that
$$\tr\B-3-\ln\det\B=\tr(e^{\Log}-\I-\Log)=\tr(\tf12\Log^2+O(\Log^3)).$$
On the other hand, we easily get
$$\tf12|\B-\I|^2=\tf12\tr(e^{2\Log}-2e^{\Log}+\I)=\tr(\tf12\Log^2+O(\Log^3)),$$
hence we also have
$$(1-\gamma)(\tr\B-3-\ln\det\B)+\tf12\gamma|\B-\I|^2=\tr(\tf12\Log^2+O(\Log^3))$$
and we see that for $\B$ being close to identity, the form of $\psi$ is almost independent of the choice of parameter $\gamma$ and the second part of $\psi$ in \eqref{helm} can be just understood as a correction for large values of $\B$.

Next, we show how the constitutive equation for $\T$ (see \eqref{T}) appears naturally if we start with the choice of the Helmholtz free energy \eqref{helm} and require that the form of the equation for $\B$ is given by \eqref{r3}. For the derivation, we followed the approach developed in \cite{Malek2018} that stems from the balance equations and requires the knowledge of how the material stores the energy, but we simplify the derivation presented there by assuming that the density is constant (in fact we set for simplicity $\varrho = 1$ and hence $\Div\ve = 0$) and the flow is isothermal, i.e., the temperature $\theta$ is constant as well. Under these assumptions the balance equations of continuum physics (for linear and angular momenta, energy and for formulation of the second law of thermodynamics) take the form
\begin{align*}
\dot{\ve}&=\Div\T, \quad \T = \T^T,\\
\dot{e}&=\T\cdot\Dv-\Div\je_e, \\
\dot{\eta}&=\xi-\Div{\je_{\eta}} \quad \textrm{ with } \xi\ge 0,
\end{align*}
where $e$ is the (specific) internal energy, $\eta$ is the entropy, $\xi$ is the  rate of entropy production, $\T$ is the Cauchy stress tensor and the quantities $\je_e$, $\je_{\eta}$ represent the internal and the entropy  fluxes, respectively. Since the quantities $\psi$, $e$, $\theta$ and $\eta$ are related through the thermodynamical identity
$$
e=\psi + \theta \eta,
$$
we can easily deduce from above identities that
\begin{equation}\label{produce}
\theta \xi= \theta\dot{\eta}+\Div{(\theta\je_{\eta})}= \T\cdot\Dv-\Div (\je_e-\theta\je_{\eta}) - \dot{\psi}.
\end{equation}
To evaluate the last term, we rewrite \eqref{r3} as
\begin{equation}\label{fuli}
	-\dot{\B}=-\la\Delta\B-a(\Dv\B+\B\Dv)-(\Wb\ve\B-\B\Wb\ve)+\delta_1(\B-\I)+\delta_2(\B^2-\B).
\end{equation}
Next, it follows from \eqref{helm} that
\begin{equation}\label{psif}
\frac{\partial \psi (\B)}{\partial \B}=\JJ,
\end{equation}
where $\JJ$ is defined as
$$\JJ:=\mu(1-\gamma)(\I-\B^{-1})+\mu\gamma(\B-\I).$$
Consequently, taking the inner product of \eqref{fuli} with $\JJ$ we observe that (since $\B \JJ = \JJ \B$, the term with $\Wb\ve$ vanishes)
\begin{equation}\label{fuli2}
\begin{split}
	-\dot{\psi}&=-\la\Delta\B\cdot \JJ-a(\Dv\B+\B\Dv)\cdot \JJ-(\Wb\ve\B-\B\Wb\ve)\cdot \JJ\\
&\qquad  +\delta_1(\B-\I)\cdot \JJ+\delta_2(\B^2-\B)\cdot \JJ\\
&=-\la \Div (\nabla \psi(\B)) -a(\Dv\B+\B\Dv)\cdot \JJ\\
&\qquad +\delta_1(\B-\I)\cdot \JJ+\delta_2(\B^2-\B)\cdot \JJ + \la \nabla \B \cdot \nabla \JJ.
\end{split}
\end{equation}
To evaluate the terms on the last line, we use the symmetry and the positive definiteness of the matrix $\B$ to obtain
\begin{equation}\label{fuli3}
\begin{split}
	(\B-\I)\cdot \JJ&=\mu (1-\gamma)|\B^{\f12}-\B^{-\f12}|^2 +\mu\gamma |\B-\I|^2,\\
(\B^2-\B)\cdot \JJ &= \mu (1-\gamma)|\B-\I|^2 +\mu \gamma |\B^{\f32}-\B^{\f12}|^2,\\
\nabla \B \cdot \nabla \JJ&=\mu \gamma |\nabla \B|^2 - \mu(1-\gamma) \nabla \B \cdot \nabla \B^{-1} \\
&= \mu \gamma |\nabla \B|^2  + \mu(1-\gamma) \nabla \B \cdot \B^{-1} \nabla \B \B^{-1}\\
&= \mu \gamma |\nabla \B|^2  + \mu(1-\gamma)|\B^{-\f12}\nn\B\B^{-\f12}|^2.
\end{split}
\end{equation}
Similarly, we obtain
\begin{equation}\label{fuli4}
a(\B\Dv + \Dv\B)\cdot \JJ = \left[ 2\mu a ((1-\gamma)(\B-\I) + \gamma (\B^2 - \B)) \right]	\cdot \Dv.
\end{equation}
Thus, using \eqref{fuli2}--\eqref{fuli4} in \eqref{produce}, we conclude that
\begin{equation}\label{xi0}
\begin{split}
\theta \xi &=-\Div (\lambda \nabla \psi (\B) + \je_e-\theta\je_{\eta})\\
&\quad +\left[ \T - 2a\mu ((1-\gamma)(\B-\I)+\gamma(\B^2-\B)) \right]\cdot \Dv \\
&\quad + \mu  \lambda (\gamma|\nabla \B|^2  + (1-\gamma) |\B^{-\f12}\nn\B\B^{-\f12}|^2)\\
 &\quad +\mu\left((1-\gamma)\delta_1|\B^{\f12}-\B^{-\f12}|^2
+\gamma\delta_2|\B^{\f32}-\B^{\f12}|^2\right)\\
&\qquad +\mu\left(((1-\gamma)\delta_2+\gamma\delta_1)|\B-\I|^2\right).
\end{split}
\end{equation}
Hence, assuming that the fluxes fulfil
\begin{equation} \label{xi0a}
\lambda \nabla \psi (\B) + \je_e-\theta\je_{\eta}=0,
\end{equation}
and setting (compare with \eqref{T})
$$
\T=-p\I+2\nu\Dv+2a\mu((1-\gamma)(\B-\I)+\gamma(\B^2-\B)),
$$
the identity \eqref{xi0} reduces to (noticing that $-p\I\cdot\Dv = -p\Div\ve = 0$)
\begin{equation}\label{xi1}
\begin{split}
\theta \xi &=\mu  \lambda (\gamma|\nabla \B|^2  + (1-\gamma) |\B^{-\f12}\nn\B\B^{-\f12}|^2)  + 2 \nu |\Dv|^2 \\
 &\quad +\mu\left((1-\gamma)\delta_1|\B^{\f12}-\B^{-\f12}|^2
+\gamma\delta_2|\B^{\f32}-\B^{\f12}|^2\right)\\
&\quad +\mu\left(((1-\gamma)\delta_2+\gamma\delta_1)|\B-\I|^2\right),
\end{split}
\end{equation}
which gives the nonnegative rate of the entropy production. Moreover, we have seen how the form of the Cauchy stress tensor $\T$ in \eqref{T} is dictated by the second line in \eqref{xi0}. Furthermore, we can also see in \eqref{xi1} (and also in the last line of \eqref{fuli3}) how the choice of the free energy \eqref{helm} affects the entropy production due to the presence of the diffusive term
$\Delta \B$ in \eqref{BB}.

\subsection{The concept of weak solution and energy (in)equality}
In order to introduce the proper concept of weak solution, we first derive the basic energy estimates based on the observations from the previous section. First, taking the scalar product of \eqref{r2} and $\ve$, we deduce the kinetic energy identity
\begin{equation}\label{kinet}
\frac12 \pp_t|\ve|^2 + \frac12 \Div(|\ve|^2 \ve) - \Div (\T \ve) + \T\cdot \D \ve = \fe \cdot\ve
\end{equation}
and replacing the term $\T\cdot \Dv$ from the equation \eqref{produce}, and using then also \eqref{xi0a} and \eqref{xi1}, we finally obtain
\begin{equation}\label{eneR}
\begin{split}
\pp_t&(\psi + \tfrac12 |\ve|^2) + \Div((\psi + \tfrac12 |\ve|^2)\ve) - \Div (\T \ve + \lambda \nabla \psi(\B) ) + 2\nu |\Dv|^2 \\
&+ \mu  \lambda \left(\gamma|\nabla \B|^2  + (1-\gamma) |\B^{-\f12}\nn\B\B^{-\f12}|^2\right)\\
 &+\mu\left((1-\gamma)\delta_1|\B^{\f12}-\B^{-\f12}|^2 +\gamma\delta_2|\B^{\f32}-\B^{\f12}|^2+((1-\gamma)\delta_2+\gamma\delta_1)|\B-\I|^2\right) \\
&\qquad = \fe \cdot\ve.
\end{split}
\end{equation}
Integrating the above identity over $\Omega$, using integration by parts and the boundary conditions \eqref{r4}--\eqref{r6}, we obtain
\begin{equation}\label{eneR2}
\begin{split}
&\f{\dd{}}{\dd{t}}\int_{\Omega} \left(\tfrac12 |\ve|^2 + \psi(\B)\right)  + 2\nu \int_{\Omega} |\Dv|^2 +\sigma \int_{\partial \Omega} |\ve |^2 \\
&\qquad +\mu  \lambda\int_{\Omega} \left( \gamma|\nabla \B|^2  + (1-\gamma) |\B^{-\f12}\nn\B\B^{-\f12}|^2\right) \\
 &\qquad\qquad +\mu\int_{\Omega}\Big((1-\gamma)\delta_1|\B^{\f12}-\B^{-\f12}|^2
+\gamma\delta_2|\B^{\f32}-\B^{\f12}|^2 \\
&\qquad\qquad\qquad\qquad\qquad\qquad\qquad +((1-\gamma)\delta_2+\gamma\delta_1)|\B-\I|^2\Big) = \int_{\Omega}\fe \cdot\ve.
\end{split}
\end{equation}
The identity \eqref{eneR2} indicates the proper choice of the function spaces for the solution $(\ve, \B)$ and the form of the (weak) formulation of the solution to \eqref{divv}--\eqref{ic2}.

\subsection{Notation}
In order to formulate the definition of a weak solution conveniently, let us fix some notation. By $L^p(\Omega)$ and $W^{n,p}(\Omega)$, $1\leq p\leq \infty$, $n\in\N$, we denote the usual Lebesgue and Sobolev space, with their usual norms denoted as $\norm{\cdot}_p$ and $\norm{\cdot}_{n,p}$, respectively. The trace operator that maps $W^{1,p}(\Omega)$ into $L^q(\partial\Omega)$, for certain $q\ge 1$, will be denoted by $\Tr$. Further, we set $W^{-1,p'}(\Omega)=(W^{1,p}(\Omega))^*$, where $p'=p/(p-1)$. We shall use the same notation for the function spaces of scalar-, vector-, or tensor-valued functions, but we will distinguish the functions themselves using different fonts such as $a$ for scalars, $\bs{a}$ for vectors and $\bb{A}$ for tensors. Also, we do not specify the meaning of the duality pairing $\langle \cdot, \cdot \rangle$, assuming that it is clear from the context. Moreover, for certain subspaces of vector valued functions, we shall use the following  notation:
\begin{align*}
C^{\infty}_{\n}&=\{\we:\Omega\to\R^3: \we\text{ infinitely differentiable},\;\we\cdot\n=0\text{ on }\partial\Omega\},\\
C^{\infty}_{\n,\Div}&=\{\we\in C^{\infty}_{\n}:\Div\we=0\text{ in }\Omega\},\\
L^2_{\n,\Div}&=\overline{C^{\infty}_{\n,\Div}}^{\norm{\cdot}_2},\quad\Wnd2=\overline{C^{\infty}_{\n,\Div}}^{\norm{\cdot}_{1,2}},\quad\Wndt=\overline{C^{\infty}_{\n,\Div}}^{\norm{\cdot}_{3,2}},\\
\Wndd2&=(\Wnd2)^*,\quad \Wnddt=(\Wndt)^*.
\end{align*}
Occasionally, we shall denote the standard inner products in $L^2(\Omega)$ and $L^2(\partial\Omega)$ as $(\cdot,\cdot)$ and $(\cdot,\cdot)_{\partial\Omega}$, respectively. The Bochner spaces of mappings from $(0,T)$ to a Banach space $X$ will be denoted as $L^p(0,T;X)$ with the norm $\norm{\cdot}_{L^p(0,T;X)}=(\int_0^T\norm{\cdot}^p_X)^{\f1p}$. If $X=L^q(\Omega)$, or $X=W^{k,q}(\Omega)$, we will write just $\norm{\cdot}_{L^pL^q}$, or $\norm{\cdot}_{L^pW^{k,q}}$, respectively. The space $\mathcal{C}_{\textrm{weak}}(0,T;X)\subset L^{\infty}(0,T;X)$ denotes a space of weakly continuous functions, i.e., for every $f\in \mathcal{C}_{\textrm{weak}}(0,T;X)$ and every $g\in X^*$ there holds
$$
\lim_{t\to t_0}\langle f(t), g\rangle = \langle f(t_0), g\rangle.
$$  The symbol $\R^{3\times 3}_{\sym}$ denotes the set of symmetric $3\times 3$ real matrices. Furthermore, by $\PD$ we denote the subset of $\R^{3\times3}_{\sym}$ which consists of positive definite matrices, i.e., those which satisfy
$$\A\z\cdot\z>0\quad\text{for all }\z\in\R^3\setminus\{0\}.$$

\section{The definition of a weak solution and its existence}\label{secmain}

In this section we state and prove the main result.

\begin{definition*}\label{DefM}
	Let $T>0$ and assume that $\Omega\subset\R^3$ is a Lipschitz domain. Let $\gamma\in(0,1)$,  $\nu,\sigma,\la>0$, $\delta_1,\delta_2\geq0$, $a\in\R$, and $\fe\in L^2(0,T;\Wndd2)$, $\ve_0\in\LH(\Omega)$. Furthermore, let $\B_0\in L^2(\Omega)$ be such that
	\begin{equation}\label{ln0}
    -\ii\ln\det\B_0<\infty.
	\end{equation}	
Then, we say that a couple $(\ve,\B):Q\to \mathbb{R}^3\times \PD$ is a weak solution to 	\eqref{divv}--\eqref{ic2} if the following hold:
\begin{align*}
\ve&\in L^2(0,T;\Wnd2)\cap L^{\infty}(0,T;L^2(\Omega)),\quad\pp_t\ve\in L^{\f43}(0,T;\Wndd2),\\
\B&\in L^2(0,T;\WH(\Omega))\cap L^{\infty}(0,T;L^2(\Omega)),\quad\pp_t\B\in L^{\f43}(0,T;\WHd(\Omega));
\end{align*}
For all $\fit\in L^{4}(0,T;\Wnd{2})$ we have
\begin{equation}
\begin{split}\label{eqv}
&\int_0^T\scal{\pp_t\ve,\fit}+\int_Q(\ve\cdot\nn)\ve\cdot\fit+\sigma\int_0^T\!\! \int_{\partial\Omega}\Tr\ve\cdot \Tr\fit\\
&\quad =-\int_Q\!(2\nu\Dv+2a\mu((1\!-\!\gamma)(\B-\I)+\gamma(\B^2-\B)))\cdot\nn\fit+\int_0^T\scal{\fe,\fit};
	\end{split}
\end{equation}
For all $\A\in L^4(0,T;\WH(\Omega))$, $\A=\A^T$, we have
\begin{equation}
    \begin{split}\label{eqB}
    &\int_0^T\scal{\pp_t\B,\A}+\int_Q((\ve\cdot\nn)\B+2\B\Wb\ve-2a\,\B\Dv)\cdot\A\\
    &\qquad\qquad +\int_Q(\delta_1(\B-\I)+\delta_2(\B^2-\B))\cdot\A+\la\int_Q\nn\B\cdot\nn\A=0;
    \end{split}
\end{equation}
The initial conditions are satisfied in the following sense
    \begin{equation}\label{ic}
    \lim_{t\to0_+}(\norm{\ve(t)-\ve_0}_2+\norm{\B(t)-\B_0}_2)=0.
    \end{equation}
    Moreover, we say that the solution satisfies the energy inequality if,  for all $t\in(0,T)$:
    \begin{equation}
    \begin{split}\label{EN}
    &\int_{\Omega} \left(\frac{|\ve(t)|^2}{2}+\psi(\B(t))\right)+\int_0^{t}\left(2\nu\norm{\Dv}_2^2+\sigma\norm{\Tr\ve}_{2,\partial\Omega}^2\right)\\
    &\qquad +\mu\lambda \int_0^t\Big((1\!-\!\gamma)\norm{\B^{-\f12}\nn\B\B^{-\f12}}_2^2+\gamma\norm{\nn\B}_2^2\Big)\\
    &\; +\mu\!\int_0^t\!\!\Big((1\!-\!\gamma)\delta_1\norm{\B^{\f12}-\B^{-\f12}}_2^2
    \!+\gamma\delta_2\norm{\B^{\f32}-\B^{\f12}}_2^2\!+ (\gamma\delta_1+(1\!-\!\gamma)\delta_2)\norm{\B-\I}_2^2\Big)\\
    &\qquad\qquad\qquad\qquad\qquad\qquad\qquad\qquad\leq \int_{\Omega}  \left(\frac{|\ve_0|^2}{2}+\psi(\B_0)\right)+\int_0^t\scal{\fe,\ve}.
    \end{split}
\end{equation}
\end{definition*}

The key result of the paper is the following
\begin{theorem*}\label{TM}
Let $T>0$ and assume that $\Omega\subset\R^3$ is a Lipschitz domain. Suppose that $\gamma\in(0,1)$,  $\nu,\sigma,\la>0$, $\delta_1,\delta_2\geq0$, $a\in\R$, and $\fe\in L^2(0,T;\Wndd2)$, $\ve_0\in\LH(\Omega)$. Furthermore, let $\B_0\in L^2(\Omega)$ be such that \eqref{ln0} holds.
Then there exists a weak solution to \eqref{divv}--\eqref{ic2} satisfying the energy inequality.
\end{theorem*}

Let us briefly explain the main difficulties connected with the analysis of the system \eqref{r1}--\eqref{r5} and our ideas how to solve them. In the standard models where $\gamma=0$, to get an a~priori estimate for $\B$, the appropriate test function to take in \eqref{r3} is $\I-\B^{-1}$. Then, using \eqref{r1} and \eqref{r2} tested by $\ve$, one can eliminate the problematic terms, such as $\B\cdot\Dv$ coming from the objective derivative. However, the non-negative quantity to be controlled, which comes from the diffusion term, turns out to be just $|\B^{-\f12}\nn\B\B^{-\f12}|^2$ and this provides little to no information. In particular, the terms $\nn\ve\B$ appearing in \eqref{r3} are going to be just integrable and it is unclear if one can show strong convergence of $\B$. Instead, one would like to test also by $\B$ to achieve control over $|\nn\B|^2$. But this is not possible, since the resulting term $\nn\ve\B\cdot\B$ cannot be estimated without some serious simplifications (such as boundedness of $\nn\ve$, two or one dimensional setting or small data). Quite remarkably, this problem is solved simply by adding $\f12\gamma |\B-\I|^2$ into the constitutive form for $\psi$. More precisely, considering  $\gamma\in(0,1)$, we observe that the appropriate test function in \eqref{r3} is in fact $(1-\gamma)(\I-\B^{-1})+\gamma(\B-\I)$. Indeed, the terms from the objective derivative cancel again due to the presence of $\gamma(\B^2-\B)$ in $\T$. But now, we also get $\gamma|\nn\B|^2$ under control, which is much better information than in the case $\gamma=0$ and it will imply compactness of all the terms appearing in \eqref{r2} and \eqref{r3}. We have seen above that such a modification of $\psi$, and consequently of $\T$, is not ad-hoc and that it rests on solid physical grounds.
	
	The second and also the last major difficulty which we will encounter is how one can justify testing of \eqref{r3} by $\B^{-1}$ on the approximate (discrete level), where $\B^{-1}$ might not even exist. This we overcome by designing a delicate approximation scheme, which takes into account the smallest eigenvalue of $\B$, and also by noting that testing \eqref{r3} only by $\B$  yields sufficiently strong a~priori estimates for the initial limit passage (in the Galerkin approximation of $\B$).
	
	Up to now, there have been no results on global existence of weak solutions to Oldroyd-B models in three dimensions, including either the standard, or diffusive variants. The closest result so far is probably \cite[Theorem~4.1]{Masmoudi2011}, however there it is assumed that $\delta_2>0$ and $\la=0$ (Giesekus model), whereas we treat also the case $\delta_2=0$, but with $\la>0$ (diffusive Oldroyd-B or Giesekus model). Moreover, in \cite{Masmoudi2011}, only the weak sequential stability of a hypothetical approximation is proved. We, on the other hand, provide the complete existence proof, including the construction of approximate solutions (which, in viscoelasticity, is generally a non-trivial task). In the article \cite{Lions2000}, Lions and Masmoudi prove the global existence in three dimensions, but only for $a=0$ (corrotational case), which is known to be much easier. The local in time existence of regular solutions for the non-diffusive variants of the models above ($\la=0$) is proved in the pioneering work \cite[Theorem~2.4.]{Saut1990}. There, also the global existence for small data is shown. In two dimensions, the problem is solved in \cite{Constantin2012} in the case $\la>0$, $\delta_1>0$, $\delta_2=0$ (diffusive Oldroyd-B model). There are also global large data existence results in three dimensions for slightly different classes of diffusive rate-type viscoelastic models, but under some simplifying assumptions. For example, in \cite{Bulicek2018} and \cite{Bulicek2019}, the authors consider the case where $\B=b\I$. This assumption, however, turns \eqref{r3} into a much simpler scalar equation. Moreover, note that if $\B=b\I$, then the equations \eqref{r2} and \eqref{r3} decouple (which is not the case in \cite{Bulicek2018} and \cite{Bulicek2019} since there the considered constitutive relation for $\T$ is more complicated than here). Furthermore, in \cite{Lukacova-Medvidova2017}, the authors consider yet another class of Peterlin viscoelastic models with stress diffusion and prove existence of a global two- or three-dimensional solution. However, the free energy associated with these models depends only on the trace of the extra stress tensor. This is a significant simplification, which can even be seen as unphysical. See also \cite{Chupin2018} for various modifications of Oldroyd-B viscoelastic models, for which an existence theory is available. Finally, in \cite{BARRETT2011} (see also \cite{Pokorny}), the global existence of a weak solution is shown for a certain regularized Oldroyd-B model (including a cut-off or nonlinear $p-$Laplace operator in the diffusive term in $\B$). Thus, one might argue that since the case $\gamma>0$ could be also seen as a regularization of the original model, we are just proving an existence of a solution to another regularization. However, this argument is not, in our opinion, correct for several reasons. First of all, the ``regularization'' $\gamma>0$ does not touch the equation \eqref{r3} at all. Second, it is not obvious why the nonlinear term $\gamma(\B-\I)^2$ should have any regularization effect. And, perhaps most importantly, we already showed in Section~\ref{shoowed} that the model with $\gamma>0$ is physically well founded and worthy of studying in its own right.
	
Since the topic is quite new and unexplored, we decided, for brevity and clarity of presentation to consider only the isothermal case. However, we believe that the framework and ideas presented here are robust enough to provide an existence analysis also for the full thermodynamical model if the evolution of the internal energy is described correctly. This is the subject of our forthcoming study.
	
\begin{remark*} Finally, we close this section with several concluding remarks on possible extensions, but we do not provide their proofs in this paper.
	\begin{itemize}
		\item[(i)] The Theorem holds also in arbitrary dimensions $d>3$ (in $d\leq2$, it is known), however with worse function spaces for the time derivatives and better for the test functions. Indeed, the only dimension-specific argument in the proof below is in the derivation of interpolation inequalities, which are then used to estimate $\pp_t\ve$ and $\pp_t\B$. Moreover, all of the non-linear terms in \eqref{eqv}, \eqref{eqB} are integrable for arbitrary $d$ if the test functions are smooth. In addition, if $d=2$, then we can prove the existence of a weak solution satisfying even the energy equality, i.e., \eqref{EN} holds with the equality sign.
		\item[(ii)] When $\Omega$ has $C^{1,1}$ boundary, then, in addition, there exists a pressure $p\in L^{\f53}(Q)$, which appears in \eqref{vv}. Then, the test functions in \eqref{eqv} need not be divergence-free if we include the term $\ii p\Div\fit$ in \eqref{eqv}. This follows in a standard way, using the Helmholtz decomposition of $\ve$ (see, e.g., \cite{Blechta2019} for details).
		\item[(iii)] It is possible to replace \eqref{r4}, \eqref{r5} by the no-slip boundary condition $\ve=\mathbf{0}$ on $\partial\Omega$. Then, we only need to change the space $\Wn{2}$ to $W^{1,2}_0$, and so on. However, then it seems that the pressure $p$ can be only obtained as a distribution (see \cite{Blechta2019}).
	\end{itemize}
\end{remark*}	

\section{Proof of the Theorem}

Throughout the proof, we shall simplify notation by assuming
$$\la=\mu=\nu=\sigma=1$$
and refer to Section~\ref{shoowed} for a detailed computation for general parameters.
To shorten all formulae, we also denote
\begin{align*}
\Sb(\A)&=(1-\gamma)(\A-\I)+\gamma(\A^2-\A) &&\textrm{for }\A\in\R^{3\times 3},\\
\R(\A)&=\delta_1(\A-\I)+\delta_2(\A^2-\A) &&\textrm{for }\A\in\R^{3\times3}.
\end{align*}

The general scheme of the proof is the following: In order to invert the matrix $\B$ and to avoid problems with low integrability in the objective derivative, we introduce the special cut-off function
$$
\rho_{\eps}(\A):=\frac{\max\{0,\La(\A)-\eps\}}{\La(\A)(1+\eps|\A|^3)} \quad \textrm{for }\A\in \R^{3\times 3}_{\sym},
$$
where $\La(\A)$ denotes a minimal eigenvalue of $\A$ (whose spectrum is real due to its symmetry)\footnote{We set $\rho_{\eps}(\A):=0$ if $\La(\A)=0$.}. Since eigenvalues of a matrix depend continuously on its entries, the function $\rho_{\eps}$ is continuous. Moreover, for any positive definite matrix $\A$ there holds $\rho_{\eps}(\A) \to 1$ as $\eps \to 0_+$. We construct a solution by an approximation scheme with parameters $k,l$ and $\eps$, where $k,l\in\N$ correspond to the Galerkin approximation for $\ve$ and $\B$, respectively, and $\eps$ corresponds to the presence of the cut-off function $\rho_{\eps}$ in certain terms. The first limit we take is $l\to \infty$, which corresponds to the limit in the equation for $\B$. This way, the limiting object $\B$ is infinite-dimensional and, using the properties of $\rho_{\eps}$, we prove that $\B^{-1}$ exists. With the help of this information, we derive the energy estimates that are uniform with respect to all the parameters. Next, we let $\eps \to 0_+$ in order to remove the truncation function and finally we take $k\to \infty$, which corresponds to the limiting procedure in the equation for the velocity $\ve$.

\subsection{Galerkin approximation} Following e.g., \cite[Appendix A.4]{Malek1996}, we know that there exists  a basis $\{\we_i\}_{i=1}^{\infty}$ of $\Wndt$, which is orthonormal in $L^2(\Omega)$ and  orthogonal in $\Wndt$. Moreover, the projection $P_k:L^2(\Omega)\to\spa\{\we_i\}_{i=1}^k$, defined as\footnote{We recall here the definition $(a,b):= \int_{\Omega} ab$.}
$$
P_k\fit=\sum_{i=1}^k(\fit,\we_i)\we_i,\quad\fit\in L^2(\Omega),
$$
is continuous in $L^2(\Omega)$ and also in $\Wndt$ independently of $k$, i.e.,
$$
\|P_k \fit\|_2\le C\|\fit\|_2 \qquad \|P_k \fit\|_{\Wndt}\le C\|\fit\|_{\Wndt}
$$
for all $\fit \in \Wndt$, where the constant $C$ is independent of $k$.  Furthermore, by the standard embedding, we also have that  $\Wndt\embl W^{2,6}(\Omega)\embl W^{1,\infty}(\Omega)$. Similarly, we construct the basis $\{\Wb_j\}_{j=1}^{\infty}$ of $\WH(\Omega)$, which is  $L^2$-orthonormal, $\WH$-orthogonal and the projection
$$
Q_l\A=\sum_{j=1}^l(\A,\Wb_j)\Wb_j,\quad\A\in L^2(\Omega),
$$
is continuous in $L^2(\Omega)$ and in $\WH(\Omega)$ independently of $l$.

Then for fixed $k,l\in\N$ and $\eps\in(0,1)$, we look for the  functions $\ve^{k,l}_{\eps},\B_{\eps}^{k,l}$ of the form
$$
\ve^{k,l}_{\eps}(t,x)=\sum_{i=1}^kc^{{k,l,\eps}}_i(t)\we_i(x)\quad\text{and}\quad \B_{\eps}^{k,l}(t,x)=\sum_{j=1}^ld^{{k,l,\eps}}_j(t)\Wb_j(x),
$$
where $c_i^{{k,l,\eps}},d_j^{{k,l,\eps}}$, $i=1,\ldots,k$, $j=1,\ldots,l$, are unknown functions of time, and we require that $\ve^{k,l}_{\eps},\B_{\eps}^{k,l}$ (and consequently the functions $c_i^{{k,l,\eps}}(t)$ and $d_j^{{k,l,\eps}}(t)$) satisfy the following system of $(k+l)$ ordinary differential equations in time interval $(0,T)$:
\begin{align}	
&\quad\begin{aligned}
\f{\dd{}}{\dd{t}}&(\ve^{k,l}_{\eps},\we_i)+((\ve^{k,l}_{\eps}\!\cdot\!\nn)\ve^{k,l}_{\eps},\we_i)+2(\D\ve^{k,l}_{\eps},\nn\we_i)+(\Tr \ve^{k,l}_{\eps},\Tr\we_i)_{\partial \Omega}\\
&\quad=-2a(\rho_{\eps}(\B^{k,l}_{\eps})\Sb(\B_{\eps}^{k,l}),\nn\we_i)+\scal{\fe,\we_i}  \qquad\textrm{for } i=1,\ldots,k,
\end{aligned}\label{gal1}\\
&\quad\begin{aligned}
\f{\dd{}}{\dd{t}}&(\B_{\eps}^{k,l},\Wb_j)+((\ve^{k,l}_{\eps}\!\cdot\!\nn)\B_{\eps}^{k,l},\Wb_j)+
(\rho_{\eps}(\B^{k,l}_{\eps})\R(\B_{\eps}^{k,l}),\Wb_j)+(\nn\B_{\eps}^{k,l},\nn\Wb_j)\\ &\quad=2(\rho_{\eps}(\B^{k,l}_{\eps})\B_{\eps}^{k,l}(a\D\ve^{k,l}_{\eps}-\Wb\ve^{k,l}_{\eps}),\Wb_j) \qquad\textrm{for } j=1,\ldots,l.
\end{aligned}
\label{galB}
\end{align}
Due to the $L^2$-orthonormality of the bases $\{\we_i\}_{i=1}^{\infty}$ and $\{\Wb_j\}_{j=1}^{\infty}$, the system \eqref{gal1}--\eqref{galB} can be rewritten as a nonlinear system of ordinary differential equations for $c^{{k,l,\eps}}_i$ and $d^{{k,l,\eps}}_j$, where $i=1,\ldots,k$ and $j=1,\ldots,l$, and we equip this system with the initial conditions
\begin{equation}\label{galic0}
c_i^{{k,l,\eps}}(0)=(\ve_0,\we_i)\qquad\text{ and }\qquad d_j^{{k,l,\eps}}(0)=(\B_0^{\eps},\Wb_j).
\end{equation}
Here, $\B_0^{\eps}$ is defined  as
$$
\B_0^{\eps}(x):=\Big\{\begin{aligned}
&\;\B_0(x)&&\text{if }\La(\B_0(x))>\eps,\\
&\;\I &&\text{elsewhere}.
\end{aligned}
$$%
Since $\B_0(x)\in \PD$ for almost every $x\in \Omega$, we have that $\Lambda(\B_0(x))>0$ for almost all $x\in \Omega$. Consequently, using the fact $\B_0\in L^2(\Omega)$, we obtain, as $\eps\to0_+$, that
\begin{equation}\label{konB0}
\begin{split}
\norm{\B_0^{\eps}-\B_0}_2^2	= \int_{\Lambda(\B_0)\leq\eps}|\I-\B_0|^2 \to 0
\end{split}
\end{equation}
Note also that the initial conditions \eqref{galic0} can be rewritten as $\ve_{\eps}^{k,l}(0)=P_k\ve_0$ and $\B_{\eps}^{k,l}(0)=Q_l\B^{\eps}_0$.

For the system \eqref{gal1}--\eqref{galic0}, Carath\'eodory's theorem can be applied and therefore there exists $T^*>0$ and absolutely continuous functions $c_i^{{k,l,\eps}}$, $d_j^{{k,l,\eps}}$  satisfying \eqref{galic0} and \eqref{gal1}--\eqref{galB} almost everywhere in $(0,T^*)$. If $T^*$ is the maximal time, for which the solution exists, and $T^*<T$, then at least one of the functions $c_i^{{k,l,\eps}}$, $d_j^{{k,l,\eps}}$ must blow up as $t\to T^*_-$. But using the estimate presented below (see \eqref{cB} valid for all $t\in(0,T^*)$), this will be seen never to happen. Thus, we can set $T^*=T$.

\subsection{Limit $l\to \infty$}\label{linfty}
In this part, we simplify the notation and denote the approximating solution, constructed in the previous section, by  $(\ve_l,\B_l):=(\ve^{k,l}_{\eps},\B^{k,l}_{\eps})$.
We start by proving estimates independent of $l$. Since $\B_l(t)$ and $\ve_l(t)$ belong for almost all $t$ to the linear hull of $\{\Wb_j\}_{j=1}^{l}$ and $\{\we_i\}_{i=1}^{k}$, respectively,  we can use $\ve_l$ instead of $\we_i$ in \eqref{gal1} and $\B_l$ instead of $\Wb_j$ in \eqref{galB} to deduce,
\begin{align}\label{est1}
\f12\f{\dd{}}{\dd{t}}&\norm{\B_l}_2^2+\norm{\nn\B_l}_2^2
=2a(\rho_{\eps}(\B_l)\B_l \Dv_l, \B_l)-(\rho_{\eps}(\B_l)\R(\B_l),\B_l),\\
\f12\f{\dd{}}{\dd{t}}&\norm{\ve_l}_2^2+2\norm{\Dv_l}_2^2+\norm{\Tr\ve_l}_{2,\partial\Omega}^2=-2a(\rho_{\eps}(\B_l)\Sb(\B_l),\Dv_l)+\scal{\fe,\ve_l},
\label{est1a}
\end{align}
where we used the integration by parts formula and the facts that $\Div \ve_l=0$ and $\Tr\ve\cdot\n=0$. Next, it follows from the definition of $\rho_{\eps}$, $\R$ and $\Sb$ that
\begin{equation}\label{hr}
\rho_{\eps}(\B_l)\left(|\Sb(\B_l)|+|\R(\B_l)||\B_l|+|\B_l|^2\right)\leq C \f{1+|\B_l|^3}{1+\eps|\B_l|^3}\leq C(\eps).
\end{equation}
Here, the notation $C(\eps)$ emphasizes that the constant $C$ depends on $\eps$; we keep this notation in what follows.
Summing \eqref{est1} and \eqref{est1a} and using the estimate \eqref{hr} to bound the term on the right-hand side, we obtain with the help of  H\"older's, Young's and Korn's inequalities that
\begin{equation}
\label{est2}	
\f{\dd{}}{\dd{t}}\left(\norm{\ve_l}_2^2+\norm{\B_l}_2^2\right)+\norm{\Dv_l}_2^2+\norm{\Tr\ve_l}_{2,\partial\Omega}^2
+\norm{\nn\B_l}_2^2\leq C(\eps)+C\norm{\fe}^2_{\Wndd2}.
\end{equation}
After integrating over $(0,T)$ with respect to time, we obtain the following bound:
\begin{equation}\label{est4}
\begin{aligned} \sup_{t\in (0,T)}\left(\norm{\ve_l}_2^2+\norm{\B_l}_2^2\right)+\int_0^{T}\left(\norm{\Dv_l}_2^2+\norm{\Tr\ve_l}_{2,\partial\Omega}^2+\norm{\nn\B_l}_2^2\right)\\\leq C(\eps)+\norm{P_k\ve_0}_2^2+\norm{Q_l\B_0^{\eps}}_2^2+C\int_0^{T}\norm{\fe}_{\Wndd2}^2\leq C(\eps),
\end{aligned}
\end{equation}
where the last inequality follows from the continuity of the projections $P_k$ and $Q_l$ and from the assumptions on data, namely that 
$$
\norm{\ve_0}_2^2+\norm{\B_0}_2^2+\norm{\ln \det \B_0}_1+C\int_0^{T}\norm{\fe}_{\Wndd2}^2<\infty.
$$

Next, we focus on the estimate for time derivatives. First, it follows from  $L^2$-orthonormality of the bases and the estimate \eqref{est4} that
\begin{equation}\label{cB}
\sum_{i=1}^kc_i(t)^2+\sum_{j=1}^ld_j(t)^2\leq C(\eps).
\end{equation}
Then, since $\ve_l$ is a linear combination of $\{\we_i\}_{i=1}^k\subset W^{1,\infty}(\Omega)$, we can estimate
\begin{equation}\label{lip}
\norm{\ve_l}_{L^{\infty}W^{1,\infty}}\leq\esssup_{t\in(0,T)}\sum_{i=1}^k|c_i(t)|\norm{\we_i}_{1,\infty}\leq C(\eps,k),
\end{equation}
and we can deduce from \eqref{gal1} that
\begin{equation}\label{dv}
\norm{\pp_t\ve_l}_{L^{\infty} W^{1,\infty}}\leq C(\eps,k).
\end{equation}
Finally, it follows from \eqref{galB} and \eqref{est4} that {\color{white}\eqref{dv}}
\begin{equation}\label{dBk}
\norm{\pp_t\B_l}_{L^{2}\WHd}\leq C(\eps,k).
\end{equation}

%
%

Using \eqref{est4}, \eqref{lip}--\eqref{dBk} and Banach-Alaoglu's theorem, we can find subsequences (which we do not relabel) and corresponding weak limits (denoted with the subscript $\eps$), such that, for $l\to\infty$, we get
\begin{align}
\ve_l&\wc\ve_{\eps} &&\text{weakly in } L^2(0,T;\Wnd2),\label{vc}\\
\ve_l&\wcs\ve_{\eps} &&\text{weakly$^*$ in }L^{\infty}(0,T;W^{1,\infty}(\Omega)),\label{wcs}\\
\pp_t\ve_l&\wcs\pp_t\ve_{\eps}&&\text{weakly$^*$ in } L^{\infty}(0,T;W^{1,\infty}(\Omega)),\label{dtc}\\
\Tr\ve_l&\wc \Tr\ve_{\eps}&&\text{weakly in } L^2(0,T;L^2(\partial\Omega)),\label{sto}\\
\B_l&\wc\B_{\eps} &&\text{weakly in } L^{2}(0,T;\WH(\Omega)),\label{lB}\\
\pp_t\B_l&\wc\pp_t\B_{\eps}&&\text{weakly in } L^{2}(0,T;\WHd(\Omega)).\label{dwB}
\end{align}
Moreover, it follows from \eqref{vc}, \eqref{dtc}, \eqref{lB}, \eqref{dwB} and from the Aubin-Lions lemma that for some further subsequences, we have{\color{white}\eqref{wcs}\eqref{sto}}
\begin{align}
\ve_l&\to\ve_{\eps} &&\text{strongly in }L^2(Q),\label{st}\\
\B_l&\to\B_{\eps} &&\text{strongly in }L^2(Q)\text{ and a.e.\ in }Q\label{stB},\\
\rho_{\eps}(\B_l)&\to \rho_{\eps}(\B_{\eps}) &&\textrm{a.e.\ in }Q.{\color{white}\eqref{st}\eqref{stB}}\label{rlp}
\end{align}
Using the convergence results \eqref{vc}--\eqref{rlp}, it is rather standard to let $l\to \infty$ in \eqref{gal1}--\eqref{galB}. This way, for almost all $t\in(0,T)$, we obtain
\begin{equation}\label{rovv}
\begin{split}
(\pp_t\ve_{\eps},\we_i)&+((\ve_{\eps}\cdot\nn)\ve_{\eps},\we_i)+2(\Dv_{\eps},\nn\we_i)+(\Tr\ve_{\eps},\Tr\we_i)_{\partial\Omega}
\\
& \quad\quad=-2a(\rho_{\eps}(\B_{\eps})\Sb(\B_{\eps}),\nn\we_i)+\scal{\fe,\we_i}
\end{split}
\end{equation}
for $i=1,\ldots,k$, and
\begin{equation}
\begin{split}
\label{rovB}
\scal{\pp_t\B_{\eps},\A} &+((\ve_{\eps}\cdot\nn)\B_{\eps},\A)+(\nn\B_{\eps},\nn\A)\\
&\quad\quad =2(\rho_{\eps}(\B_{\eps})\B_{\eps}(a\Dv_{\eps}-\Wb\ve_{\eps}),\A)-(\rho_{\eps}(\B_{\eps})\R(\B_{\eps}),\A)
\end{split}
\end{equation}
for all $\A\in\WH(\Omega)$. Moreover, from \eqref{lB} and \eqref{dwB}, we get $\B_{\eps}\in \mathcal{C}(0,T; L^2(\Omega))$ and it is standard to show that $\B_{\eps}(0,\cdot)=\B_0^{\eps}$ and $\ve_{\eps}(0,\cdot)=P_k\ve_0$.

\subsection{Limit $\eps\to 0$} In this part we consider the solutions $(\ve_{\eps}, \B_{\eps})$ constructed in the preceding section for $\eps \in (0,1)$ and we study their behaviour as $\eps \to 0_+$. To do so, we first have to derive estimates that are uniform with respect to $\eps$. Following the ideas used before in the derivation of the model, we wish to test \eqref{rovB} by the function
\begin{equation}\label{testa}
\JJ_{\eps}:=(1-\gamma)(\I-\B_{\eps}^{-1})+\gamma(\B_{\eps}-\I).
\end{equation}
This test function, however, contains $\B^{-1}_{\eps}$ and we need to justify that it exists (for any $\eps\in(0,1)$).

\subsubsection{Estimates for the inverse matrix - still $\eps$-dependent}
First, we prove that $\La(\B_{\eps})\ge \eps$. For this purpose, let $\z\in\R^3$ be arbitrary and consider\footnote{In this subsection, we use the notation $(f)_+:=\max\{0,f\}$ and $(f)_-:=\min\{0,f\}$.}
\begin{equation}\label{test}
\A=(\B_{\eps}\z\cdot\z-\eps|\z|^2)_{-}\,(\z\ot\z),\quad\text{where }(\z\ot\z)_{ij}:=z_iz_j
\end{equation}
in \eqref{rovB}. Due to the properties of $\B_{\eps}$ (see \eqref{lB}), we know that $\A$ belongs to $L^2(0,T;\WH(\Omega))$ and we can use it as a test function in \eqref{rovB}. Upon inserting $\A$ into \eqref{rovB}, we integrate the result over $(0,\tau)$ with some fixed $\tau\in(0,T)$. We evaluate all terms in \eqref{rovB} separately. For the time derivative, we have
\begin{equation}\label{du1}
\begin{split}
\int_0^{\tau} \scal{\pp_t\B_{\eps},\A}&=\int_0^{\tau}\scal{\pp_t(\B_{\eps}\z\cdot\z-\eps|\z|^2), (\B_{\eps}\z\cdot\z-\eps|\z|^2)_{-}}\\
&=\f12\norm{(\B_{\eps}(\tau)\z\cdot\z-\eps|\z|^2)_-}_2^2-\f12\norm{(\B^{\eps}_0\z\cdot\z-\eps|\z|^2)_-}_2^2\\
&=\f12\norm{(\B_{\eps}(\tau)\z\cdot\z-\eps|\z|^2)_-}_2^2,
\end{split}
\end{equation}
where, for the last equality, the definition of $\B^{\eps}_0$ was used.
Furthermore, we obtain
\begin{equation}\label{du2}
\begin{split}
\int_Q\nn\B_{\eps}\cdot\nn\A&=\int_0^{\tau}\ii\nn(\B_{\eps}-\eps\I)\cdot\nn((\B_{\eps}\z\cdot\z-\eps|\z|^2)_-\,(\z\ot\z))\\
&=\int_0^{\tau}\norm{\nn(\B_{\eps}\z\cdot\z-\eps|\z|^2)_-}_2^2
\end{split}
\end{equation}
and
\begin{equation}\label{du3}
\begin{split}
\int_Q(\ve_{\eps}\cdot\nn)\B_{\eps}\cdot\A&=\int_0^{\tau}\ii\ve_{\eps}\cdot\nn(\B_{\eps}\z\cdot\z-\eps|\z|^2)(\B_{\eps}\z\cdot\z-\eps|\z|^2)_-\\
&=\frac12\int_0^{\tau}\ii\ve_{\eps}\cdot\nn((\B_{\eps}\z\cdot\z-\eps|\z|^2)_{-}^2\\
&=-\frac12\int_0^{\tau}\ii ((\B_{\eps}\z\cdot\z-\eps|\z|^2)_{-}^2\Div \ve_{\eps}=0,
\end{split}
\end{equation}
integrating by parts and using the fact that $\Div \ve_{\eps}=0$ and $\Tr\ve_{\eps}=0$.
Since
$$\B_{\eps}\z\cdot\z\geq \La(\B_{\eps})|\z|^2\qquad \textrm{a.e. in }Q,$$
we also observe, that
\begin{multline}
0\geq(\La(\B_{\eps})-\eps)_+(\B_{\eps}\z\cdot\z-\eps|\z|^2)_{-}
\geq(\La(\B_{\eps})-\eps)_+(\La(\B_{\eps})-\eps)_{-}\,|\z|^2=0.
\end{multline}
Hence, we get
\begin{equation}
\rho_{\eps}(\B_{\eps})\A=0 \qquad \textrm{a.e.\ in }Q.\label{du0}
\end{equation}
Consequently, inserting $\A$ of the form \eqref{test} into \eqref{rovB}, we see that the right-hand side is identically zero. Therefore, relations \eqref{du1}, \eqref{du2}, \eqref{du3} and \eqref{du0} yield
$$
\begin{aligned}
&\norm{(\B_{\eps}\z\cdot\z-\eps|\z|^2)_-}_2^2(\tau)\\
&\qquad \le \norm{(\B_{\eps}\z\cdot\z-\eps|\z|^2)_-}_2^2(\tau)+2\int_0^{\tau}\norm{\nn(\B_{\eps}\z\cdot\z-\eps|\z|^2)_-}^2_{2}=0,
\end{aligned}
$$
which implies
\begin{equation}\label{pze}
\B_{\eps}\z\cdot\z\geq\eps|\z|^2\quad\text{for every }\z\in\R^3\text{ and a.e.\ in }Q.
\end{equation}
Thus, we have the following estimate for the minimal eigenvalue of $\B_{\eps}$:
$$\La(\B_{\eps})\geq\inf_{0\neq\z\in\R^3}\f{\B_{\eps}\z\cdot\z}{|\z|^2}\geq\eps.$$
Therefore, the inverse matrix $\B_{\eps}^{-1}$ is well defined and satisfies
\begin{equation}\label{det}
|\B_{\eps}^{-1}|\le \frac{C}{\eps}\quad\text{a.e.\ in Q.}
\end{equation}
Furthermore, since
$$
\begin{aligned}
\nabla \B_{\eps}^{-1}&=\B_{\eps}^{-1}\B_{\eps}\nabla \B_{\eps}^{-1}=\B_{\eps}^{-1}\nabla (\B_{\eps} \B_{\eps}^{-1})-\B_{\eps}^{-1} (\nabla \B_{\eps}) \B_{\eps}^{-1}=-\B_{\eps}^{-1} (\nabla \B_{\eps}) \B_{\eps}^{-1},
\end{aligned}
$$
we conclude from \eqref{est4} and \eqref{det}, that
\begin{equation}
\int_Q |\nabla \B_{\eps}^{-1}|^2 \le \int_Q |\B_{\eps}^{-1}|^4 |\nabla \B_{\eps}|^2 \le C(\eps). \label{inverse}
\end{equation}
Hence, the inverse of $\B_{\eps}$ exists and $\B^{-1}_{\eps}\in L^2(0,T;\WH(\Omega))$.

\subsubsection{Estimates independent of $(\eps,k)$}
At this point, we can test \eqref{rovB} with $\JJ_{\eps}$ defined in \eqref{testa}. This way, we obtain
\begin{equation}
\begin{split}
\label{rovBJ}
&\scal{\pp_t\B_{\eps},\JJ_{\eps}} +((\ve_{\eps}\cdot\nn)\B_{\eps},\JJ_{\eps})+(\nn\B_{\eps},\nn\JJ_{\eps})\\
&\qquad =2(\rho_{\eps}(\B_{\eps})\B_{\eps}(a\Dv_{\eps}-\Wb\ve_{\eps}),\JJ_{\eps})-(\rho_{\eps}(\B_{\eps})\R(\B_{\eps}),\JJ_{\eps}).
\end{split}
\end{equation}
Next, we evaluate all terms. Here, we follow very closely the procedure developed in Section~\ref{shoowed}, see the derivation of \eqref{fuli2} and consequent identities. Since
$$
\JJ_{\eps} = \frac{\partial \psi(\B_{\eps})}{\partial \B_{\eps}},
$$
where $\psi$ is defined in \eqref{helm}, it is clear that
\begin{align}
\label{time12}
\scal{\pp_t\B_{\eps},\JJ_{\eps}}&=\frac{\dd{}}{\dd{t}}\ii \psi(\B_{\eps}),\\
((\ve_{\eps}\cdot\nn)\B_{\eps},\JJ_{\eps})&= \ii \ve_{\eps} \cdot \nabla \psi(\B_{\eps})=0.\label{time13}
\end{align}
Next, recalling \eqref{fuli3}, we get
\begin{equation}\label{fuli3B}
\begin{aligned}
(\rho_{\eps}(\B_{\eps})\R(\B_{\eps}),\JJ_{\eps})&=\ii \rho_{\eps}(\B_{\eps})\left(\delta_1(1-\gamma)|\B_{\eps}^{\f12}-\B_{\eps}^{-\f12}|^2 +(\delta_1\gamma +\delta_2(1-\gamma))|\B_{\eps}-\I|^2 \right.\\
&\qquad \left.+\delta_2\gamma |\B_{\eps}^{\f32}-\B_{\eps}^{\f12}|^2\right),\\
(\nabla \B_{\eps} ,\nabla \JJ_{\eps})
&= \gamma \|\nabla \B_{\eps}\|_2^2  + (1-\gamma)\|\B_{\eps}^{-\f12}\nn\B_{\eps}\B_{\eps}^{-\f12}\|_2
^2
\end{aligned}
\end{equation}
and due to the fact that $\B_{\eps}\JJ_{\eps}=\JJ_{\eps}\B_{\eps}$ we also have
$$
\begin{aligned}
(\rho_{\eps}(\B_{\eps})(\Wb\ve_{\eps}\B_{\eps}-\B_{\eps}\Wb\ve_{\eps}),\JJ_{\eps})&=0,\\
a(\rho_{\eps}(\B_{\eps})(\Dv_{\eps}\B_{\eps}+\B_{\eps}\Dv_{\eps}),\JJ_{\eps})&=2a (\rho_{\eps}(\B_{\eps})\D \ve_{\eps}, \B_{\eps}\JJ_{\eps})\\
&=2a (\rho_{\eps}(\B_{\eps})\D \ve_{\eps}, (1-\gamma)(\B_{\eps}-\I)+\gamma(\B_{\eps}^2-\B_{\eps}))\\
&=2a(\rho_{\eps}(\B_{\eps})\Sb(\B_{\eps}), \D \ve_{\eps}),
\end{aligned}
$$
where we used the fact that the trace of $\D\ve_{\eps}$ is identically zero.
Hence, using $\A:=\JJ_{\eps}$ (defined in \eqref{testa}) in \eqref{rovB} and taking into account the above identities, we deduce that
\begin{equation}\label{ESThlp}
\begin{split}
\frac{\dd{}}{\dd{t}}\ii \psi(\B_{\eps})+(1-\gamma)&\norm{\B_{\eps}^{-\f12}\nn\B_{\eps}\B^{-\f12}_{\eps}}_2^2+\gamma\norm{\nn\B_{\eps}}_2^2\\
+(\gamma\delta_1+(1-\gamma)\delta_2)&\norm{\sqrt{\rho_{\eps}(\B_{\eps})}(\B_{\eps}-\I)}_2^2\\
+(1-\gamma)\delta_1&\norm{\sqrt{\rho_{\eps}(\B_{\eps})}(\B_{\eps}^{\f12}-\B_{\eps}^{-\f12})}_2^2\\
+\gamma\delta_2&\norm{\sqrt{\rho_{\eps}(\B_{\eps})}(\B_{\eps}^{\f32}-\B_{\eps}^{\f12})}_2^2=2a(\rho_{\eps}(\B_{\eps})\Sb(\B_{\eps}), \D \ve_{\eps}).
\end{split}
\end{equation}

Similarly as in previous section, replacing $\we_i$ in \eqref{rovv} by $\ve_{\eps}$, we get
\begin{align}
\f12\f{\dd{}}{\dd{t}}\norm{\ve_{\eps}}_2^2+2\norm{\Dv_{\eps}}_2^2+\norm{\Tr\ve_{\eps}}_{2,\partial\Omega}^2&=\scal{\fe,\ve_{\eps}}-2a(\rho_{\eps}(\B_{\eps})\Sb(\B_{\eps}),\Dv_{\eps}).
\label{est1aa}
\end{align}
Thus, summing \eqref{ESThlp} and \eqref{est1aa} and integrating the result with respect to time $t\in (0,\tau)$, we deduce the identity
\begin{equation}\label{EST}
\begin{split}
&\f12\norm{\ve_{\eps}(\tau)}^2_2+\ii\psi(\B_{\eps}(\tau))\\
&\qquad+\int_0^{\tau}\Big(2\norm{\Dv_{\eps}}_2^2+\norm{\Tr\ve_{\eps}}_{2,\partial\Omega}^2
+(1-\gamma)\norm{\B_{\eps}^{-\f12}\nn\B_{\eps}\B^{-\f12}_{\eps}}_2^2+\gamma\norm{\nn\B_{\eps}}_2^2\\
&\qquad +(\gamma\delta_1+(1-\gamma)\delta_2)\norm{\sqrt{\rho_{\eps}(\B_{\eps})}(\B_{\eps}-\I)}_2^2\\
&\qquad+(1-\gamma)\delta_1\norm{\sqrt{\rho_{\eps}(\B_{\eps})}(\B_{\eps}^{\f12}-\B_{\eps}^{-\f12})}_2^2
+\gamma\delta_2\norm{\sqrt{\rho_{\eps}(\B_{\eps})}(\B_{\eps}^{\f32}-\B_{\eps}^{\f12})}_2^2\Big)
\\
&=\f12\norm{P_k\ve_0}^2_2+\ii\psi(\B^{\eps}_0)+\int_0^{\tau}\scal{\fe,\ve_{\eps}}\leq \f12\norm{\ve_0}^2_2+\ii\psi(\B_0)+\int_0^{\tau}\scal{\fe,\ve_{\eps}},
\end{split}
\end{equation}
where, for the last inequality we used the continuity of $P_k$, the definition of $\B_{0}^{\eps}$ and the fact that $\psi(\I)=0$.

From \eqref{EST}, we get, using Korn's, Sobolev's, H\"older's and Young's inequalities, that
\begin{equation}
\label{aprioriew}
\norm{\ve_{\eps}}_{L^{\infty}L^2}+\norm{\ve_{\eps}}_{L^2L^6}+\norm{\ve_{\eps}}_{L^2W^{1,2}}+\norm{\B_{\eps}}_{L^2W^{1,2}}+\norm{\B_{\eps}}_{L^2L^6}
\leq C,
\end{equation}
where the constant $C$ depends only on $\Omega$, $\ve_0$, $\B_0$ and $\fe$. Furthermore, the interpolation inequalities yield
\begin{equation}\label{inter}
\norm{\ve_{\eps}}_{L^{\f{10}3}L^{\f{10}3}}+\norm{\ve_{\eps}}_{L^4L^3}+\norm{\B_{\eps}}_{L^{\f{10}3}L^{\f{10}3}}+\norm{\B_{\eps}}_{L^4L^3}
+\norm{\B_{\eps}}_{L^{\f83}L^4}\leq C.
\end{equation}

Finally, we focus on the estimate for time derivatives.
Let $\fit\in L^4(0,T;\Wndt)$ be such that $\norm{\fit}_{L^4W^{3,2}}\leq 1$. Then, since $\ve_{\eps}$ is a linear combination of $\{\we_i\}_{i=1}^k$, we obtain, using \eqref{rovv}, H\"older's inequality, \eqref{EST}, \eqref{inter} and $W^{3,2}$-continuity of $P_k$, that
$$\int_0^T\scal{\pp_t\ve_{\eps},\fit}\leq C,$$
hence
\begin{equation}\label{Dv}
\norm{\pp_t\ve_{\eps}}_{L^{\f43}\Wnddt}\leq C.
\end{equation}
Similarly, by considering $\A\in L^4(0,T;\WH(\Omega))$ in \eqref{rovB}, we get
\begin{equation}\label{DB}
\norm{\pp_t\B_{\eps}}_{L^{\f43}\WHd}\leq C.
\end{equation}

\subsubsection{Limit $\eps\to0_+$.} From \eqref{aprioriew}, \eqref{Dv}, \eqref{DB}, Banach-Alaoglu's theorem and the Aubin-Lions lemma, we obtain the existence of a couple $(\ve_k, \B_k)$ satisfying the following convergence results\footnote{The convergence results \eqref{strv}, \eqref{strB} are true in any space $L^p(Q)$, $1\leq p<\f{10}3$, as can be seen from \eqref{inter} and Vitali's theorem. The space $L^3(\Omega)$ is chosen for simplicity; in our proof, we need $p>2$.}
\begin{align}
\ve_{\eps}&\wc\ve_k && \text{weakly in}\quad L^2(0,T;\Wnd2),\\
\pp_t\ve_{\eps}&\wc\pp_t\ve_k&&\text{weakly in}\quad L^{\f43}(0,T;\Wnddt),\\
\Tr\ve_{\eps}&\wc\Tr\ve_k&&\text{weakly in}\quad L^2(0,T;L^2(\partial\Omega)),\\
\B_{\eps}&\wc\B_k && \text{weakly in}\quad L^2(0,T;\WH(\Omega)),\\
\pp_t\B_{\eps}&\wc\pp_t\B_k&&\text{weakly in}\quad L^{\f43}(0,T;\WHd(\Omega)),\\
\ve_{\eps}&\to\ve_k &&\text{strongly in }L^3(Q)\text{ and a.e.\ in }Q\label{strv},\\
\B_{\eps}&\to\B_k &&\text{strongly in }L^3(Q)\text{ and a.e.\ in }Q\label{strB}.\\
\end{align}
Using \eqref{strB} and letting $\eps\to0_+$ in \eqref{pze}, we obtain
\begin{equation}\label{pzk}
\B_k\z\cdot\z\geq0\quad\text{a.e.\ in }Q\text{ and for all }\z\in\R^3.
\end{equation}
Hence $\La(\B_k)\geq0$ and $\det\B_k\geq0$ a.e.\ in $Q$. Therefore, using \eqref{strB} again and the continuity of $\psi$, there exists (still possibly infinite) limit
\begin{align*}
\psi(\B_{\eps})&\to\psi(\B_k)\quad\text{a.e.\ in }Q.
\end{align*}
However, since $\psi\geq0$, Fatou's lemma implies that, for almost every $t\in(0,T)$, we have
$$\ii\psi(\B_k)(t)\leq\liminf_{\eps\to0_+}\ii\psi(\B_{\eps})(t)\leq C.$$
Thus, we deduce that
\begin{equation}\label{lndetn}
\norm{\psi(\B_{k})}_{L^{\infty}L^1}\leq C.
\end{equation}
If there existed a set $E\subset Q$ of a positive measure, where $\La(\B_k)=0$, then also $-\ln\det\B_k=\infty$ on that set, which contradicts \eqref{lndetn}.
Thus, we have
\begin{equation}\label{lam}
\La(\B_k)>0\text{ a.e.\ in }Q.
\end{equation}
Therefore, it directly follows from the continuity of $\La$, that $\rho_{\eps}(\B_{\eps})\to1$ a.e.\ in $Q$. Then, since $\rho_{\eps}(\B_{\eps})\leq1$, we further get, by Vitali's theorem, that
\begin{align}
\rho_{\eps}(\B_{\eps})&\to 1\quad\text{strongly in }L^p(Q) \textrm{ for all } p\in [1,\infty).\label{convre}
\end{align}
Using the established convergence results, it is easy to let $\eps\to0_+$ in \eqref{rovv} and \eqref{rovB} and obtain, for almost all $t\in(0,T)$, that
\begin{equation}
\begin{split}\label{rovv2}
&\scal{\pp_t\ve_k,\we_i}+((\ve_k\cdot\nn)\ve_k,\we_i)+2(\Dv_k,\nn\we_i)\\
&\qquad =-(\Tr\ve_k,\Tr\we_i)_{\partial\Omega}-2a(\Sb(\B_k),\nn\we_i)+\scal{\fe,\we_i},\quad \text{for }i=1,\ldots,k,
\end{split}
\end{equation}
and that
\begin{equation}
\begin{split}\label{rovB2}
&\scal{\pp_t\B_k,\A}+((\ve_k\cdot\nn)\B_k,\A)+(\nn\B_k, \nn\A)\\
&\qquad =2(\B_k(a\Dv_k-\Wb\ve_k),\A)-(\R(\B_k),\A)\quad\text{for all }\A\in\WH(\Omega).
\end{split}
\end{equation}
Furthermore, we can take the limit in the estimates \eqref{EST}, \eqref{inter}, \eqref{Dv} and \eqref{DB} using either the weak lower semi-continuity of norms or, in the terms which depend on $\B_{\eps}$, e.g.\ $\int_Q\rho_{\eps}(\B_{\eps})|\B_{\eps}^{\f32}-\B_{\eps}^{\f12}|^2$, we apply \eqref{lam} to conclude the pointwise limit and then use Fatou's lemma. Thus, inequalities \eqref{EST}, \eqref{inter}, \eqref{Dv} and \eqref{DB} continue to hold in the same form, but for $(\ve_k,\B_k)$ instead of $(\ve_{\eps},\B_{\eps})$ and with $1$ instead of $\rho_{\eps}(\B_{\eps})$. In particular, for almost all $t\in(0,T)$, we have
\begin{equation}\label{ESTkk}
\begin{split}
&\f12\norm{\ve_{k}(\tau)}^2_2+\ii\psi(\B_{k}(\tau))\\
&\quad+\int_0^{\tau}\Big(2\norm{\Dv_{k}}_2^2+\norm{\Tr\ve_{k}}_{2,\partial\Omega}^2
+(1-\gamma)\norm{\B_{k}^{-\f12}\nn\B_{k}\B^{-\f12}_{k}}_2^2+\gamma\norm{\nn\B_{k}}_2^2\\
&\qquad +(\gamma\delta_1+(1-\gamma)\delta_2)\norm{\B_{k}-\I}_2^2\\
&\qquad +(1-\gamma)\delta_1\norm{\B_{k}^{\f12}-\B_{k}^{-\f12}}_2^2
+\gamma\delta_2\norm{\B_{k}^{\f32}-\B_{k}^{\f12}}_2^2\Big)
\\
&\leq \f12\norm{\ve_0}^2_2+\ii\psi(\B_0)+\int_0^{\tau}\scal{\fe,\ve_k}.
\end{split}
\end{equation}
The attainment of initial conditions is standard (see the last section for details in a more complicated case).

\subsection{Limit $k\to\infty$}
Since we start from the same a~priori estimates as in the previous section, we follow, step by step, the procedure developed when taking the limit $\eps\to0_+$. The only difference is that the term $\rho_{\eps}(\B_{\eps})$ is not present. Thus, using the density of $\{\we_i\}_{i=1}^{\infty}$ in $\Wndt$, we obtain, after letting $k\to\infty$, for almost all $t\in(0,T)$, that
\begin{equation}
\begin{split}\label{rovv3}
&\scal{\pp_t\ve,\fit}+((\ve\cdot\nn)\ve,\fit)+2(\Dv,\nn\fit)\\
&\qquad =-(\Tr\ve,\Tr\fit)_{\partial\Omega}-2a(\Sb(\B),\nn\fit)+\scal{\fe,\fit}\quad\text{for all }\fit\in\Wndt
\end{split}
\end{equation}
and that
\begin{equation}
\begin{split}\label{rovB3}
&\scal{\pp_t\B,\A}+((\ve\cdot\nn)\B,\A)+(\nn\B, \nn\A) \\
&\qquad=2(\B(a\Dv-\Wb\ve),\A)-(\R(\B),\A)\quad\text{for all }\A\in\WH(\Omega).
\end{split}
\end{equation}

Moreover, from the weak lower semi-continuity of norms, we obtain the energy inequality~\eqref{EN} for almost all $t\in (0,T)$. Furthermore, the same argument as above implies that $\B$ is positive definite a.e.\ in $Q$. Now observe that, by H\"older's inequality and \eqref{inter}, all the terms in \eqref{rovv3} except the first one, are integrable for every $\fit\in L^4(0,T;\Wnd2)\embl L^4(0,T;L^6(\Omega))$. Indeed, for example for the non-linear terms, we get
$$\int_Q|(\ve\cdot\nn)\ve\cdot\fit|\leq\norm{\ve}_{L^4L^3}\norm{\nn\ve}_{L^2L^2}\norm{\fit}_{L^4L^6}$$
and
$$\int_Q|\Sb(\B)\cdot\nn\fit|\leq C\norm{\B}_{L^{\f83}L^4}^2\norm{\nn\fit}_{L^4L^2}.$$
Hence, the functional $\pp_t\ve$ can be uniquely extended to $\pp_t\ve\in L^{\f43}(0,T;\Wndd2)$ and we can use the density argument to conclude \eqref{eqv}. Analogously, we obtain \eqref{eqB}. Hence, it remains to show that \eqref{EN} holds for all $t\in (0,T)$ and that the initial data fulfil \eqref{ic}.

\subsubsection{Energy inequality for all $t\in (0,T)$}
First, we observe, that due to \eqref{aprioriew}, \eqref{Dv} and \eqref{DB}, we have that
\begin{equation}\label{weakc}
\begin{aligned}
\ve &\in \mathcal{C}_{\textrm{weak}}(0,T; L^2(\Omega))\quad\text{and}\quad
\B \in \mathcal{C}_{\textrm{weak}}(0,T; L^2(\Omega)).
\end{aligned}
\end{equation}
Next, we notice that the function $\psi$ is convex on the convex set $\PD$. Indeed, evaluating the second Fr\'echet derivative of $\psi$, we get
$$
\frac{\partial^2\psi(\A)}{\A^2}=(1-\gamma)\A^{-1}\ot\A^{-1}+\gamma\I\ot\I\quad\text{for all }\A\in\PD,$$
which is obviously a positive definite operator for any $\gamma\in[0,1]$ and consequently, $\psi$ must be convex on $\PD$.

Further, we integrate \eqref{EN} over $(t_1,t_1+\delta)$, where $t_1\in(0,T)$, and divide the result by $\delta$. Using also an elementary inequality
$$\int_0^{t_1}g\leq\f1{\delta}\int_{t_1}^{t_1+\delta}\left(\int_0^tg\right)\dd{t}$$
valid for every integrable non-negative $g$, we get
\begin{equation}\label{ESTint}
\begin{split}
&\frac{1}{2\delta}\int_{t_1}^{t_1+\delta}\norm{\ve(t)}^2_2+\frac{1}{\delta}\int_{t_1}^{t_1+\delta}\ii\psi(\B(t))\\
&\qquad+\int_0^{t_1}\Big(2\norm{\Dv}_2^2+\norm{\Tr\ve}_{2,\partial\Omega}^2
+(1-\gamma)\norm{\B^{-\f12}\nn\B\B^{-\f12}}_2^2+\gamma\norm{\nn\B}_2^2\\
&\qquad\qquad\qquad\qquad +(\gamma\delta_1+(1-\gamma)\delta_2)\norm{\B-\I}_2^2\\
&\qquad\qquad\qquad\qquad\qquad +(1-\gamma)\delta_1\norm{\B^{\f12}-\B^{-\f12}}_2^2
+\gamma\delta_2\norm{\B^{\f32}-\B^{\f12}}_2^2\Big)
\\
&\qquad\qquad\qquad\qquad\qquad\qquad\qquad\leq \f12\norm{\ve_0}^2_2+\ii\psi(\B_0)+\frac{1}{\delta}\int_{t_1}^{t_1+\delta}\int_0^{\tau}\scal{\fe,\ve}.
\end{split}
\end{equation}
Finally, we let $\delta\to 0_+$. The limit on the right hand side is standard and consequently, if we show that
\begin{equation}\label{needs}
\begin{split}
\frac{1}{2}\norm{\ve(t_1)}^2_2+\ii\psi(\B(t_1))\le \liminf_{\delta\to 0_+} \frac{1}{\delta}\int_{t_1}^{t_1+\delta}\left(\frac{\norm{\ve(t)}^2_2}{2}+\ii\psi(\B(t))\right),
\end{split}
\end{equation}
then \eqref{EN} will hold for all $t\in (0,T)$. To show it, we notice that due to \eqref{weakc}
\begin{equation}\label{contqw}
\begin{aligned}
\ve(t) &\rightharpoonup  \ve(t_1) &&\textrm{weakly in $L^2(\Omega)$ as $t\to t_1$},\\
\B(t) &\rightharpoonup  \B(t_1) &&\textrm{weakly in $L^2(\Omega)$ as $t\to t_1$},\\
\end{aligned}
\end{equation}
Consequently, due to the weak lower semicontinuity and the convexity of $\psi$ we also have for all $t\in (0,T)$
$$
\ii |\ve(t)|^2 + \psi(\B(t)) \le C.
$$
Hence denoting by $\Omega_M\subset \Omega$ the set where $|\ve(t_1,\cdot)|+|\B(t_1,\cdot)| + |\B^{-1}(t_1,\cdot)| \le M$, it follows from the previous estimate that $|\Omega \setminus \Omega_M| \to 0$ as $M\to \infty$. Hence, since $\psi$ is nonnegative and convex, we have for all $t\in (t_1,t_1+\delta)$ that
$$
\begin{aligned}
&\ii \frac{|\ve(t)|^2}{2} +\psi(\B(t)) \ge \int_{\Omega_M}\!\! \frac{|\ve(t)|^2}{2} +\psi(\B(t))\\
&\ge \int_{\Omega_M}\!\!\!\! \frac{|\ve(t_1)|^2}{2} +\psi(\B(t_1))+ \int_{\Omega_M}\!\! \!\! \ve(t_1) \cdot (\ve(t)-\ve(t_1)) + \frac{\partial \psi(\B(t_1))}{\partial \B} \cdot (\B(t)-\B(t_1)).
\end{aligned}
$$
Since, $\ve(t_1)$ and $\partial_{\B}\psi(\B(t_1))$ are bounded on $\Omega_M$, we can integrate the above estimate over $(t_1,t_1+\delta)$ and it follows from \eqref{contqw} that
$$
\begin{aligned}
&\liminf_{\delta \to 0_+} \frac{1}{\delta}\int_{t_1}^{t_1+\delta}\ii \frac{|\ve(t)|^2}{2} +\psi(\B(t)) \ge \int_{\Omega_M} \frac{|\ve(t_1)|^2}{2} +\psi(\B(t_1)).
\end{aligned}
$$
Hence, letting $M\to \infty$, we deduce \eqref{needs} and the proof of \eqref{EN} is complete.

\subsubsection{Attainment of initial conditions}
First, it is standard to show from the construction and from the weak continuity \eqref{contqw}, that for arbitrary $\fit,\A\in L^2(\Omega)$ there holds
\begin{equation}\label{v0}
\begin{aligned}
\lim_{t\to0_+}(\ve(t),\fit)&=(\ve_0,\fit)\quad\text{and}\quad\lim_{t\to0_+}(\B(t),\A)=(\B_0,\A).
\end{aligned}
\end{equation}
Next, using the convexity of $\psi$ and \eqref{v0} (and consequently weak lower semicontinuity of the corresponding integral) and letting $t\to 0_+$ in \eqref{EN}, we deduce that
\begin{equation}\label{con-norm}
\begin{split}
\norm{\ve_0}_2^2&+2\ii\psi(\B_0)\le \liminf_{t\to0_+}\left(\norm{\ve(t)}_2^2+2\ii\psi(\B(t))\right) \\
&\le \limsup_{t\to0_+}\left(\norm{\ve(t)}_2^2+2\ii\psi(\B(t))\right)\le \norm{\ve_0}_2^2+2\ii\psi(\B_0).
\end{split}
\end{equation}
We claim that this implies that
\begin{equation}\label{con-norm2}
\begin{split}
\norm{\ve_0}_2^2=\lim_{t\to0_+} \norm{\ve(t)}_2^2\quad\text{and}\quad\ii\psi(\B_0)&=\lim_{t\to0_+} \ii\psi(\B(t)).
\end{split}
\end{equation}
Indeed, assume for a moment that
$$
\norm{\ve_0}_2^2 < \liminf_{t\to0_+}\norm{\ve(t)}_2^2.
$$
But then it follows from \eqref{con-norm} that
$$
\ii \psi(\B_0) > \liminf_{t\to0_+} \ii \psi(\B(t)),
$$
which contradicts \eqref{v0} and convexity of $\psi$. Consequently, \eqref{con-norm2} holds.

It directly follows from \eqref{v0}$_1$ and \eqref{con-norm2}$_1$ that
$$
\lim_{t\to 0_+}\norm{\ve(t)-\ve_0}_2^2=0.
$$
To claim the same result also for $\B$, we simply split $\psi$ as follows
$$
\psi(\A)=\frac{\gamma}{2} |\A - \I|^2 + (1-\gamma) (\tr \A - 3 - \ln \det \A)=:\gamma \psi_1(\A) + (1-\gamma) \psi_2(\A).
$$
Similarly as above, it is easy to observe that $\psi_1$ as well as $\psi_2$ are convex on the set of positive definite matrices. Therefore, \eqref{con-norm2}$_2$ and \eqref{v0}$_2$ imply
\begin{equation}\label{con-norm3}
\begin{split}
\ii |\B_0-\I|^2&=2\ii \psi_1(\B_0)=2\lim_{t\to 0_+} \ii \psi_1(\B(t))= \lim_{t\to 0_+} \ii |\B(t)-\I|^2,\\
\ii\psi_2(\B_0)&=\lim_{t\to0_+} \ii\psi_2(\B(t)).
\end{split}
\end{equation}
Finally, \eqref{v0} and \eqref{con-norm3}$_1$ lead to
$$
\begin{aligned}
\lim_{t\to 0_+}\norm{\B(t)-\B_0}_2^2&=\lim_{t\to 0_+}\norm{(\B(t)-\I) +(\I-\B_0)}_2^2\\
&=\lim_{t\to 0_+}\left(\norm{\B(t)-\I}_2^2 +\norm{\B_0-\I}_2^2-2\ii (\B(t)-\I)\cdot (\B_0 -\I)\right)\\
&=0,
\end{aligned}
$$
which finishes the proof of \eqref{ic} and consequently also the proof of the Theorem.


\begin{thebibliography}{10}
	\bibitem{Amirat}
	 Y. Amirat, D. Bresch, J. Lemoine, J. Simon, Effect of rugosity on a flow governed by stationary
	Navier-Stokes equations, \emph{Quart. Appl. Math.} \textbf{59} (2001), no.~4, 769--785.
	
	\bibitem{Amrouche2014}
	C. Amrouche, A. Rejaiba, ${\rm L}^p$-theory for Stokes and Navier–Stokes equations with Navier boundary condition, \emph{J.
	Differential Equations} \textbf{256} (2014), no.~4, 1515--1547.

	\bibitem{Amrouche2011}
	C. Amrouche, N.E.H. Seloula, On the Stokes equations with the Navier-type boundary conditions, \emph{Differ. Equ. Appl.} \textbf{3}
	(2011), no.~4, 581--607.
	
	\bibitem{Baba2019}
	H. Al Baba, {M}aximal $L^p$-$L^q$ regularity to the Stokes problem with Navier boundary conditions, 
	\emph{Adv. Nonlinear Anal.} \textbf{8} (2019), no.~1, 743--761. 
	
	\bibitem{BARRETT2011}
	J. Barrett, S. Boyaval, {E}xistence and approximation of a (regularized)
		{O}ldroyd-{B} model, \emph{Math. Models Methods Appl. Sci.}
	\textbf{21} (2011), no.~09, 1783--1837.
	
	\bibitem{Basson}
	A. Basson, D. G\'erard-Varet
	Wall laws for fluid flows at a boundary with random roughness,
	\emph{Comm. Pure Appl. Math.} \textbf{61} (2008), no.~7, 941--987. 
	
	\bibitem{Beirao2020}
	H. Beir\~{a}o da Veiga, J. Yang, {R}egularity criteria for {N}avier-{S}tokes equations with slip
	boundary conditions on non-flat boundaries via two velocity
	components, \emph{Adv. Nonlinear Anal.} 
	\textbf{9} (2020), no.~1, 633--643.
	
	\bibitem{Blechta2019}
	J. Blechta, J. M\'{a}lek, K.R. Rajagopal, {O}n the Classification of Incompressible Fluids and a~Mathematical Analysis of the Equations That Govern Their Motion,
	\emph{SIAM J. Math. Anal.} \textbf{52} (2020), no.~2, 1232--1289. 
	
	\bibitem{Bucur}
	D. Bucur, E. Feireisl,
	The incompressible limit of the full Navier-Stokes-Fourier system on domains with rough boundaries,
	\emph{Nonlinear Anal. Real World Appl.} \textbf{10} (2009), no.~5, 3203--3229. 
	
	\bibitem{Bulicek2019}
	M. Bul{\'{\i}}{\v{c}}ek, E. Feireisl, J. M{\'{a}}lek, On a class of
		compressible viscoelastic rate-type fluids with stress-diffusion,
	\emph{Nonlinearity} \textbf{32} (2019), no.~12, 4665--4681.
	
	\bibitem{Bulicek2009}
	M. Bul{\'{\i}}{\v{c}}ek, J. M{\'{a}}lek, K.R. Rajagopal, Mathematical
		analysis of unsteady flows of fluids with pressure, shear-rate, and
		temperature dependent material moduli that slip at solid boundaries, \emph{{SIAM} J. Math. Anal.} \textbf{41} (2009), no.~2, 665--707.
	
	\bibitem{stick1}
	M. Bul\'{\i}\v{c}ek, J. M\'{a}lek, Internal flows of incompressible
		fluids subject to stick-slip boundary conditions, \emph{Vietnam J. Math.}
	\textbf{45} (2017), no.~1-2, 207--220.
	
	\bibitem{stick2}
	\bysame, Large data analysis for {K}olmogorov's two-equation model of
		turbulence, \emph{Nonlinear Anal. Real World Appl.} \textbf{50} (2019), 104--143.
	
	\bibitem{Bulicek2018}
	M. Bul\'{\i}\v{c}ek, J. M\'{a}lek, V. Pr{\r{u}}\v{s}a, E. S\"{u}li, P{DE}
		analysis of a class of thermodynamically compatible viscoelastic rate-type
		fluids with stress-diffusion, \emph{Contemp. Math.}, vol. 710, Amer. Math.
	Soc., Providence, RI, 2018, pp.~25--51.
	
	\bibitem{Zaba3}
	M. Bul\'{\i}\v{c}ek, J. M\'{a}lek, J. \v{Z}abensk\'{y}, On generalized
		{S}tokes' and {B}rinkman's equations with a pressure-and shear-dependent
		viscosity and drag coefficient, \emph{Nonlinear Anal. Real World Appl.} \textbf{26}
	(2015), 109--132.
	
	\bibitem{Zaba2}
	M. Bul\'{\i}\v{c}ek, J. \v{Z}abensk\'{y}, Large data existence theory for
		unsteady flows of fluids with pressure- and shear-dependent viscosities,
	\emph{Nonlinear Anal.} \textbf{127} (2015), 94--127.
	
	\bibitem{Chupin2018}
	L. Chupin, Global strong solutions for some differential viscoelastic
		models, \emph{{SIAM} J. Appl. Math.} \textbf{78} (2018), no.~6,
	2919--2949.
	
	\bibitem{Constantin2012}
	P. Constantin, M. Kliegl, Note on global regularity for two-dimensional
		{O}ldroyd-{B} fluids with diffusive stress, \emph{Arch. Ration. Mech. Anal.} \textbf{206} (2012), no.~3, 725--740.
	
	\bibitem{Dostalik2019}
	M. Dostal\'\i{}k, V. Pr{\r{u}}\v{s}a, T. Sk\v{r}ivan, On diffusive
		variants of some classical viscoelastic rate-type models, \emph{AIP Conference
	Proceedings} \textbf{2107} (2019).

	\bibitem{Drda}
	S.-Q. Wang, P.A. Drda, Stick–slip transition in capillary flow of polyethylene. 2. Molecular weight dependence and low temperature
	anomaly, \emph{Macromolecules} \textbf{29} (11) (1996) 4115--4119.
	
	\bibitem{Elkareh}
	A.W. El-Kareh, L.G. Leal, Existence of solutions for all {D}eborah
		numbers for a non-{N}ewtonian model modified to include diffusion, \emph{J. Non-Newton. Fluid Mech.} \textbf{33} (1989), no.~3, 257--287.
	
	\bibitem{Saut1990}
	C. Guillop\'{e}, J.C. Saut, Existence results for the flow of
		viscoelastic fluids with a differential constitutive law, \emph{Nonlinear Anal.} \textbf{15} (1990), no.~9, 849--869.
		
	\bibitem{Hatzi}
	S.G. Hatzikiriakos, Wall slip of molten polymers, \emph{Prog. Polym. Sci.} \textbf{37} (2012) 624--643.
	
	\bibitem{Hron2017}
	J. Hron, V. Milo{\v{s}}, V. Pr{\r{u}}{\v{s}}a, O. Sou{\v{c}}ek, K. T{\r{u}}ma,
	On thermodynamics of incompressible viscoelastic rate type fluids with
		temperature dependent material coefficients, \emph{Internat. J. Non-Linear Mech.} \textbf{95} (2017), 193--208.
		
	\bibitem{Jager}
	W. J\"ager, A. Mikeli\'c,
	On the roughness-induced effective boundary conditions for an incompressible viscous flow,
	\emph{J. Differential Equations} \textbf{170} (2001), no.~1, 96--122. 
	
	\bibitem{Pokorny}
	O. Kreml, M. Pokorn\'{y}, P. \v{S}alom, On the global existence for a
		regularized model of viscoelastic non-{N}ewtonian fluid, \emph{Colloq. Math.} \textbf{139} (2015), no.~2, 149--163.
	
	\bibitem{Leray}	
	J. Leray, Sur le mouvement d'un liquide visqueux emplissant l'espace,
	\emph{Acta Math.} \textbf{63} (1934), no.~1, 193--248. 
	
	\bibitem{Lions2000}
	P.L. Lions, N. Masmoudi, Global solutions for some oldroyd models of
		non-newtonian flows, \emph{Chin. Ann. Math. Ser. B} \textbf{21} (2000),
	no.~2, 131--146.
	
	\bibitem{Lukacova-Medvidova2017}
	M. Luk{\'{a}}{\v{c}}ov{\'{a}}-Medvi{\v{d}}ov{\'{a}}, H. Mizerov{\'{a}},
	{\v{S}}. Ne{\v{c}}asov{\'{a}}, M. Renardy, Global existence result for
		the generalized {P}eterlin viscoelastic model, \emph{{SIAM} J. Math. Anal.} \textbf{49} (2017), no.~4, 2950--2964.
		
	\bibitem{Macha2014}
	V. M\'acha, J. Tich\'y,
	Higher integrability of solutions to generalized Stokes system under perfect slip boundary conditions, \emph{J. Math. Fluid Mech.} \textbf{16} (2014), no.~4, 823--845.
	
	\bibitem{Malek1996}
	J. M{\'{a}}lek, J. Ne{\v{c}}as, M. Rokyta, M. R{\r{u}}{\v{z}}i{\v{c}}ka,
	\emph{{W}eak and {M}easure-valued {S}olutions to {E}volutionary {PDE}s},
	Chapman \& Hall, 1996.
	
	\bibitem{Malek2018}
	J. M{\'{a}}lek, V. Pr{\r{u}}{\v{s}}a, T. Sk{\v{r}}ivan, E. S\"{u}li,
	Thermodynamics of viscoelastic rate-type fluids with stress diffusion,
	\emph{Phys. Fluids} \textbf{30} (2018).
	
	\bibitem{MaPr18}
	J. M\'{a}lek, V. Pr{\r{u}}\v{s}a, \emph{Derivation of Equations for Continuum
		Mechanics and Thermodynamics of Fluids}, Handbook of mathematical analysis in
	mechanics of viscous fluids, Springer, Cham, 2018, pp.~3--72.
	
	\bibitem{MaRaTu15}
	J. M\'{a}lek, K.R. Rajagopal, K. T{\r{u}}ma, On a variant of the
		{M}axwell and {O}ldroyd-{B} models within the context of a thermodynamic
		basis, \emph{Internat. J. Non-Linear Mech.} \textbf{76} (2015), 42--47.
	
	\bibitem{MaRaTu18}
	\bysame, Derivation of the variants of the {B}urgers model using a
		thermodynamic approach and appealing to the concept of evolving natural
		configurations, \emph{Fluids} \textbf{3} (2018), no.~4.
	
	\bibitem{Zaba4}
	E. Maringov\'{a}, J. \v{Z}abensk\'{y}, On a
		{N}avier-{S}tokes-{F}ourier-like system capturing transitions between viscous
		and inviscid fluid regimes and between no-slip and perfect-slip boundary
		conditions, \emph{Nonlinear Anal. Real World Appl.} \textbf{41} (2018), 152--178.
	
	\bibitem{Masmoudi2011}
	N. Masmoudi, Global existence of weak solutions to macroscopic models of
		polymeric flows, \emph{J. Math. Pures Appl. (9)}
	\textbf{96} (2011), no.~5, 502--520.
	
	\bibitem{RS_2000}
	K.R. Rajagopal, A.R. Srinivasa, A thermodynamic frame work for rate type
		fluid models, \emph{J. Non-Newton. Fluid Mech.} \textbf{88} (2000), no.~3,
	207--227.
	
	\bibitem{RS_2004}
	K.R. Rajagopal, A.R. Srinivasa, On thermomechanical restrictions of
		continua, \emph{Proc. R. Soc. Lond., Ser. A, Math. Phys. Eng. Sci.} \textbf{460}
	(2004), no.~2042, 631--651.
	
\end{thebibliography}

\providecommand{\bysame}{\leavevmode\hbox to3em{\hrulefill}\thinspace}
\providecommand{\MR}{\relax\ifhmode\unskip\space\fi MR }
\providecommand{\MRhref}[2]{%
	\href{http://www.ams.org/mathscinet-getitem?mr=#1}{#2}
}
\providecommand{\href}[2]{#2}

\end{document}